\documentclass[12pt]{book}
\usepackage{jfaa2e}     
\usepackage{latexsym}
\usepackage{amsmath}
\usepackage{amsthm}
\usepackage{amssymb}
\usepackage{amscd}
\usepackage[english]{babel}
\usepackage{graphicx,amssymb}
\usepackage[usenames]{color}

\font\msbm=msbm10
\newcommand\vect[1]{{\bf#1}}
\newcommand\matr[1]{{\bf#1}}

\def \mathbb#1{\hbox{\msbm{#1}}}

\newtheorem{thm}{Theorem}[section]
\newtheorem*{thmA}{Theorem A}
\newtheorem*{thmB}{Conjecture B}
\newtheorem*{thmC}{Theorem C}

\newtheorem{lem}[thm]{Lemma}
\newtheorem{prop}[thm]{Proposition}
\newtheorem{defn}[thm]{Definition}

\numberwithin{equation}{section}

\newcommand{\mt}{\mathbb}   
\newcommand{\mc}{\mathcal}  

\newcommand{\ra}    { \rightarrow }
\newcommand{\ds}    {\displaystyle}
\def        \half   {{\frac{1}{2}}}

\doi{jfaatest- 05/28/04} 

\begin{document}

\author{Joseph Shtok and Michael Elad}

\chapter{Analysis of Basis Pursuit Via Capacity Sets}

\begin{center}
The Computer Science Department, \\The Technion - Israel Institute
of Technology,  Haifa 32000 Israel, \\ email:
[shtok,elad]@cs.technion.ac.il.
\end{center}

\footnotetext{\textit{Math Subject Classifications:} 68P30, 68W25.}

\footnotetext{\textit{Keywords and Phrases:} sparse representations,
$\ell_1$-reconstruction, Basis Pursuit, random support, capacity
sets.}

\begin{abstract}
Finding the sparsest solution $\vect{\alpha}$ for an under-determined
linear system of equations $\matr{D}\vect{\alpha}=\vect{s}$ is of
interest in many applications. This problem is known to be NP-hard.
Recent work studied conditions on the support size of $\vect{\alpha}$
that allow its recovery using $\ell_1$-minimization, via the Basis
Pursuit algorithm. These conditions are often relying on a scalar
property of $\matr{D}$ called the mutual-coherence. In this work we
introduce an alternative set of features of an arbitrarily given
$\matr{D}$, called the {\bf capacity sets}. We show how those could
be used to analyze the performance of the basis pursuit, leading to
improved bounds and predictions of performance. Both theoretical and
numerical methods are presented, all using the capacity values, and
shown to lead to improved assessments of the basis pursuit success in
finding the sparest solution of $\matr{D}\vect{\alpha}=\vect{s}$.
\end{abstract}


\section{Introduction}\label{intro}
A powerful trend in signal processing that has evolved in recent
years is the use of redundant dictionaries, rather than just bases,
for a sparse representation of signals (images, sound tracks, and
more). In such a setting, we consider a linear equation
$\vect{s}=\matr{D}\vect{\alpha}$, where $\vect{s}$ is a given signal,
$\matr{D}$ is the representation dictionary, and $\vect{\alpha}$ is
the signal's representation. The matrix $\matr{D}$ is a general full
rank $N\times L$ matrix, where $L>N$, assumed to have $\ell_2$
normalized columns. The number of non-zero elements in the
coefficient vector $\vect{\alpha}$  is measured by the $\ell_0$-norm,
$\|\cdot\|_0$, on $\mt{R}^L$. The goal is to find, within the
$(L-N)$-dimensional affine space of the solutions for this equation,
the sparsest representation for $\vect{s}$, i.e. one which has the
least number of non-zero entries. This goal is formalized by the
following optimization problem:
$$ \ds
(P_0):\ \ \ \ \ \mbox{Arg}\min_{\vect{\alpha} \in
\mt{R}^L}\|\vect{\alpha}\|_0\ \ s.t. \ \
\matr{D}\vect{\alpha}=\vect{s} .
$$
In this paper, we consider the signals for which the solution of
$(P_0)$ is unique, and we define $\mc{S}(\matr{D})$ as the family of
such signals. We denote $\Omega=\{1,...,L\}$, and refer to the
support of the vector $\vect{\alpha}=(\alpha_1,...,\alpha_L)^T$ as
the set $\Gamma = supp(\vect{\alpha})=\{n\in \Omega\ |\ \alpha_n\neq
0 \}$.

The problem $(P_0)$ is NP-hard, demanding an exhaustive search over
all the subsets of columns of $\matr{D}$ \cite{NP-Hard}. One of the
most effective techniques to approximate its solution is the convex
relaxation of the $\ell_0$-norm. It uses the $\ell_1$-norm, the
closest convex norm on $\mt{R}^L$:
$$\ds
(P_1):\ \ \ \ \ \mbox{Arg} \min_{\vect{\alpha}\in
\mt{R}^L}\|\vect{\alpha}\|_1\ \ s.t.\ \
\matr{D}\vect{\alpha}=\vect{s}.
$$
The solution of $(P_1)$ is carried out by linear programming. We are
interested in signals $\vect{s}\in \mc{S}(\matr{D})$ for which the
solutions of $(P_0)$ and $(P_1)$ coincide. The idea of using $(P_1)$
to find the sparsest solution is called Basis Pursuit (BP), as
coined by Chen, Donoho and Saunders ~\cite{Ch,CDS}.

Let $\vect{\alpha}$ be a representation of $\vect{s}$, with support
$\Gamma=supp(\vect{\alpha})\subset\Omega$. The matrix
$\matr{D}_{\Gamma}$ is a matrix of size $N\times |\Gamma|$
containing the columns (also referred to as atoms) of $\matr{D}$
used for the construction of $\vect{s}$. This matrix is necessarily
full-rank (with rank equals $|\Gamma|$). Knowing the support
$\Gamma$ suffices to enable perfect recovery of $\vect{\alpha}$, and
thus our interest is confined to the ability to recover the support
$\Gamma$.

\begin{defn}
A subset $\Gamma\subset\Omega$ is called
\textbf{$\ell_1$-reconstructible} with respect to the dictionary
$\matr{D}$ if the solution of $(P_1)$ coincides with the solution of
$(P_0)$ for every signal $s\in \mc{S}(\matr{D})$ that admits a
representation with the support $\Gamma$.
\end{defn}

The main task of the paper is to obtain  conditions on support sizes
which imply that they are  $\ell_1$-reconstructible. For any specific
support $\Gamma\subset\Omega$ there exists a straightforward (yet
exhaustive) test whether it admits recovery by BP -- simply apply BP
to the finite family of signals $\vect{s}=\matr{D}\vect{\alpha}$
generated from coefficient vectors $\vect{\alpha}$ with the support
$\Gamma$ covering all possible sign patterns (i.e. $2^{|\Gamma|}$
such tests\footnote{In fact, half of this amount is required because
if $\vect{\alpha}$ is reconstructible, then so is
$-\vect{\alpha}$.}). If the recovery succeeds for all these choices
of $\vect{\alpha}$, it will also succeed for any other representation
with support $\Gamma$ \cite{DonohoHuo,Malioutov}.

Clearly, such a testing approach is impractical in most cases. If we
aim to find the prospects of success of the BP for a fixed
cardinality $|\Gamma|$, this requires a set of tests as described
above per each possible support $\Gamma$ having such a cardinality,
and this implies a need for approximately $L^{|\Gamma|}$ groups of
tests. Thus, the exhaustive approach should be replaced  either by a
random set of tests with empirical claims, or a theoretical study.

Within the theoretical  attempts to estimate the power of the BP, two
approaches are distinguished in the existing literature. Earlier work
carried out the worst case analysis for a given dictionary, providing
conditions on the support cardinality that guarantee that any support
satisfying them is $\ell_1$-reconstructible
\cite{DE1,DonohoHuo,EB1,Fuchs1,Gribonval_UNION,T03}. These conditions
are often very restrictive and far from empirical evidence. Another,
more recent, approach presents a probabilistic analysis, providing
conditions for special families of dictionaries under which
\emph{most} signals of a given cardinality are
$\ell_1$-reconstructible \cite{CRT1,CR2,D4,DT1,Tropp06a}. The results
depict a general asymptotic behavior with regard to the sparse
support recovery.

In both worst-case and probabilistic-analysis branches of work, many
classical results rely heavily on a scalar feature of the dictionary,
known as the {\em mutual-coherence}
\cite{DE1,Fuchs1,Gribonval_UNION,T03}. A related measure also used is
the Babel function \cite{DE1,T03}. More recent work employs the
Restricted Isometry Property (RIP) \cite{RIP}. The information
carried by all these measures is very pessimistic; furthermore, the
RIP is very expensive computationally and mainly used for theoretical
analysis. In this work we set to improve the existing worst case
results for a given general dictionary $\matr{D}$, as reported in
\cite{DE1,Fuchs1,Gribonval_UNION,T03}. We achieve this progress by
replacing the above-mentioned with a set of alternative features that
we refer to as the {\em capacity sets} of the dictionary. A thorough
computational analysis of $\matr{D}$ and probabilistic tools are
applied to the problem, leading to improved probabilistic bounds.

In the next section we recall the existing theoretical results
concerning $\ell_1$-recovery as a function of the support
cardinality. In section 3 we define two versions of the {\em
capacity set} and present the main theoretical results of this paper
using these features. Section 4 expands on the above results by
providing two numerical algorithms using the {\em capacity sets}.
Section 5 provides an overall comparison of the various methods
presented in this work to assess the performance of BP for several
test-cases.


\section{Background}\label{background}

Most known results on sparsity rely on the \emph{mutual-coherence},
denoted as $\mu$, of the dictionary. This is the maximum of the inner
products between the columns: $\mu=\max_{i\neq
j\in\Omega}|<\vect{d}_i,\vect{d}_j>|$. This correlation between the
columns, reflected in its worst value by $\mu$, helps establishing
the "safe zone" for the support sizes, where both the uniqueness of
sparsest representation and its $\ell_1$-recovery can be guaranteed.

For $\matr{D}=[\matr{\Phi}_1,\matr{\Phi}_2]$ a pair of orthonormal
bases, the following sufficient condition for $\Gamma$ to be
$\ell_1$-reconstructible is proven in \cite{EB1}:
$$
|\Gamma|\leq \frac{\sqrt{2}-0.5}{\mu}~~.
$$
Donoho and Elad in ~\cite{DE1} treat a general dictionary
$\matr{D}$. They define the problem
\begin{equation}\label{C_gamma}\ds
(C_{\Gamma}):\ \ \ \ \ \max_{\vect{\delta}\in Null(\matr{D})}
\sum_{k\in \Gamma}|\delta_k|\ \ \ s.t. \ \ \|\vect{\delta}\|_1=1~,
\end{equation}
and show that its solution is intimately tied to the ability to
recover the support $\Gamma$, by the following lemma:
\begin{lem}\label{lem_DE1}(\cite{DE1}, \textit{Lemma 2})
A sufficient condition on the support $\Gamma$ to be
$\ell_1$-reconstructible is
\begin{equation}
val(C_{\Gamma})<\half.
\end{equation}
\end{lem}
\noindent This criteria is used to prove the following theorem:
\begin{thm}\label{thm_DE1}(\cite{DE1}, \textit{Theorem 7})
A sufficient condition on a support $\Gamma\subset\Omega$ to be
$\ell_1$-reconstructible is
\begin{equation}\label{CB}
|\Gamma|<\half\left(1+\frac{1}{\mu}\right).
\end{equation}
\end{thm}

Typically, the coherence behaves at best like
$\mc{O}(\frac{1}{\sqrt{N}})$, hence the results stated above predict
quite weak $\ell_1$-recovery, which is refuted by the empirical
evidence: usually BP recovers supports of size proportional to $N$
(and not its squared-root).

A generalization of the coherence is introduced in \cite{DE1} and
later used by J. Tropp in ~\cite{T03}: for any $0\le m \le L$, the
Babel function $\mu_1(m)$ is defined by
$$\ds
\mu_1(m)=\max_{|\Lambda|=m}~~\max_{\eta\in\Omega\backslash
\Lambda}\sum_{\lambda\in\Lambda}
|<\vect{\phi}_{\lambda},\vect{\phi}_{\eta}>|.
$$
In terms of this function, a support of size $m$ is proven to be
$\ell_1$-reconstructible provided the following inequality holds
\cite{T03}:
$$
\mu_1(m-1)+\mu_1(m)<1.
$$
Unfortunately, in cases where the coherence $\mu$ is close to $1$
(implying an existence of at least one problematic pair of atoms),
the growth of $\mu_1(m)$ is too fast to provide any improvement.

Average case analysis improves the asymptotic bounds on
reconstructible support sizes. The work in \cite{CR2} shows that for
the dictionary $\matr{D}=[\matr{I},\matr{F}^*]$, where $\matr{F}$ is
the Fourier transform, random uniformly sampled support admits
$\ell_1$-recovery with high probability if (the expectation of) its
cardinality is $\mc{O}(N/\log N)$, which improves the
$\mc{O}(\sqrt{N})$ estimation of the worst case approach. For a
general orthonormal pair, it is shown in (\cite{CR2}, Theorem 5.3)
that most random supports which cardinality behaving like
$\mc{O}(1/(\mu^2 \log ^6 N))$ admit recovery by BP. The $\log N$
appearing in these expressions is suspected by the authors of
\cite{CR2} to be unnecessary, which in effect turns this expression
into $\mc{O}(N)$ (for incoherent dictionaries). A similar and related
result, exhibiting the square of the mutual coherence in the
denominator of the bound, appears in \cite{Tropp06a}. As such, this
result is effective in cases where the dictionary is ``uniformly
coherent'', and the methods employed are not very suitable for
dictionaries with high coherence.

The idea that representations with cardinalities $\mc{O}(N)$ are
$\ell_1$ -reconstructible is supported by the results reported in
\cite{D4,Donoho04a,DT1}. This result is obtained for asymptotically
growing dictionaries of size $N\times \delta N$ constructed by
concatenating random vectors of unit $l_2$-norm, independently drawn
from the uniform distribution. It is shown that all supports of size
up to $\rho(\delta) N$ are $\ell_1$-reconstructible with probability
approaching 1. The work in \cite{Donoho04a,DT1} provides theoretical
assessments for $\rho(\delta)$, based on connection to study on
neighborly polytopes. Despite being asymptotical, these results
illuminate the empirically-supported evidence regarding the
reconstruction abilities of minimal $L_0$-norm  supports by linear
programming.

As good as these results sound, they do not provide useful numerical
information about the ability of $\ell_1$-reconstruction applied to a
specifically given dictionary $\matr{D}$ of certain size, which is a
practical and central question in the application of BP. Such
information can only be obtained today by results involving the
coherence $\mu$ or its descendants. Thus, the gap is especially big
when the dictionary is not uniformly coherent and when $\mu\gg
\dfrac{1}{\sqrt{N}}$.

In this work we introduce new features of the dictionary $\matr{D}$,
the {\em capacity sets}. These features are obtained as the solutions
to specific linear programming problems that probe the dictionary
$\matr{D}$. We consider two such options: a vector of capacities
$\vect{q}$ and a matrix $\matr{Q}$, as we shall explain in details in
the next section. These features are used to develop novel analysis
of BP performance as a function of the support's cardinality.

One interesting benefit of the proposed analysis is a better
treatment of dictionaries which are not ``uniformly coherent''. In
cases where there exists a small set of columns in $\matr{D}$ with
strong linear dependency, the coherence and the babel function behave
badly, tending to lead to overly pessimistic bounds. As we show, the
use of the capacities leads in these cases to much better results.
Besides that, the capacities are shown to be more delicate indicators
of the dictionary, as reflected in a better prediction of the BP
performance.

Use of {\em capacity sets} bridges the gap between purely theoretical
estimations of the reconstructible support sizes for given dictionary
$\matr{D}$, which are usually fast but provide pessimistic lower
bound, and the empirical tests of  $\matr{D}$, which give very
accurate account on BP-reconstruction abilities, but are
computationally prohibitive. We propose theoretical results and
algorithms that employ the {\em capacity sets} to perform
computational assessment of these abilities, which is fast relative
to full empirical test and more optimistic than known practical
formulae. The question of computational complexity is discussed in
details in section \ref{sect_complexity}.



\section{Capacity Sets and Their Use}\label{overview}

In this section we define two versions of the {\em capacity sets},
and state the main theoretical results that employ them for the
analysis of the BP.

\subsection{The Capacity Vector $\vect{q}$}

The capacity vector consists of elements related to  an intermediate
tool used in the proof of Theorem \ref{thm_DE1} in \cite{DE1}:

\begin{defn}\label{CV}
The capacity vector $\vect{q}=(q_1,...,q_L)^T$ of a dictionary $\matr{D}\in
\mt{R}^{N\times L}$ is defined for all $k\in\Omega$ by
\begin{equation}\label{defn_cvector}
 q_k=\max_{\delta\in Null(D)}\delta_k \ \ \ s.t. \ \
\|\delta\|_1=1.
\end{equation}
\end{defn}

Computing the elements of $\vect{q}$ is relatively easy, and amounts
to a simple set of $L$ independent linear programming problems of the
form
$$
\hat{\vect{x}}_k = \mbox{Arg}\min_{\vect{x}}~||\vect{x}||_1
~~~\mbox{subject to}~~~\matr{D}\vect{x}=\vect{0} ~~\mbox{and}~~ x_k
= 1,$$ and then assigning $q_k = 1/||\hat{\vect{x}}_k||_1$.

To see the equivalence of the two problems, notice that the vector
$\tilde{\vect{x}}_k =\hat{\vect{x}}_k/\|\hat{\vect{x}}_k\|_1$ is an
element of null space of $\matr{D}$ with unit $\ell_1$-norm. Since
$(\hat{\vect{x}}_k)_k=1$ and $\|\hat{\vect{x}}_k\|_1$ is smallest
possible, the value $q_k
=1/||\hat{\vect{x}}_k||_1=(\tilde{\vect{x}}_k)_k$ is just the
solution of \ref{defn_cvector}.


Via Lemma \ref{lem_DE1}, the definition of $\vect{q}$ provides a
sufficient condition $\sum_{k\in \Gamma}q_k<\half$ on a given support
$\Gamma$ to ensure its recovery by $\ell_1$-minimization.
Furthermore, by gathering the $|\Gamma|$ largest entries from
$\vect{q}$, a simple generalization of Theorem 2.2 can be proposed.
However, in this work we seek a better bound that takes into account
the variety of possible supports, rather than the worst one. One such
numerical technique is suggested in section 4, proposing a special
quantization of the values in $\vect{q}$ to obtain a lower bound on
the fraction of support sizes which admit recovery by BP.

In this section we aim to obtain a more theoretically flavored
result that uses $\vect{q}$. Denote by $E_q$ the mean value of the
capacity vector $\vect{q}$, and by $\sigma^2_q$ its variance
$\frac{1}{L}\sum_{k\in \Omega}(q_k-E_q)^2$. The following theorem
uses these quantities to evaluate the probability of
$\ell_1$-reconstruction for a given support size:

\begin{thmA}
For any $1\leq \ell < \frac{1}{2E_q}$, a support $\Gamma$ of size
$\ell$, sampled uniformly at random from $\Omega$, admits
$\ell_1$-recovery with probability
\begin{eqnarray}\label{eq:BoundA}
P(\ell) > \frac{  \left( \half -\ell E_q\right)^2}{{\ell
\sigma_q^2}+\left( \half -\ell E_q\right)^2}~.
\end{eqnarray}
\end{thmA}
In the special case of a constant capacity vector, the theorem boils
down to support size threshold of $\frac{1}{2E_q}$, since then the
variance becomes zero. We show in Section
\ref{sect_back_to_coherence} that weakened version of Theorem A
yields the classical threshold of
$|\Gamma|<\half\left(1+\frac{1}{\mu}\right)$ (see Theorem
\ref{thm_DE1}).

\noindent {\bf Proof}: We fix $\ell$ and chose subsets $\Lambda,
\Gamma \subset \Omega$ according to two different probability models.
The elements of $\Gamma$ are chosen uniformly from $\Omega$ without
replacement and form a set of $\ell$ distinct column indices. The
$\ell$ elements of $\Lambda$ are chosen uniformly with replacement
(i.e. $\Lambda$ is a multiset of size $\ell$ with possible
duplicates). Now, define random variables
\begin{equation}\label{defn_xl_yl}
x_{\ell}=\sum_{k\in \Gamma}q_k,\;\;\; y_{\ell}=\sum_{m\in \Lambda}q_m.
\end{equation}
In these terms, the probability $P(\ell)$, defined in the statement
of the theorem, is bounded below by $P(x_{\ell}< \half)$. In turn,
we shall bound the probability $P(x_{\ell}< \half)$ by means of the
Tchebychev inequality, which involves the mean and the variance of
$x_{\ell}$. These parameters are easily computable for $y_{\ell}$:
by its definition, we have $\mt{E}(y_{\ell})=\ell E_q$,
$var(y_{\ell})=\ell\sigma_q^2$. Our result is based on the following
connection between the variables $x_{\ell}$ and $y_{\ell}$, as shown
in Appendix A:
\begin{equation}\label{X_Y_rel1}
\mt{E}(x_{\ell})=\mt{E}(y_{\ell})~~\mbox{and}~~ var(x_{\ell})\leq
var(y_{\ell}).
\end{equation}
Given any real scalar $a>0$, the one-tailed version of the
Tchebychev inequality \cite{Polya} for $x_{\ell}$ reads
$$
P\left(x_{\ell}-E_x\geq a\sigma_x\right) = P\left(x_{\ell}\geq E_x+
a\sigma_x\right)\leq \frac{1}{1+a^2}~,
$$
where $E_x=\mt{E}(x_{\ell})$, $\sigma^2_x=var(x_{\ell})$.

By (\ref{X_Y_rel1}), we substitute $E_x=\ell E_q$. Also, since a
larger variance implies a lower probability, we put
$\sqrt{\ell}\sigma_q$ instead of $\sigma_x$ and obtain
$$
P\left(x_{\ell}\geq \ell E_q+ a\sqrt{\ell}\sigma_q\right)\leq
P\left(x_{\ell}\geq E_x+ a\sigma_x\right)\leq \frac{1}{1+a^2}.
$$

The parameter $a$ is chosen such that $\ell E_q+
a\sqrt{\ell}\sigma_q = \half$, leading to $a=(\half -\ell
E_q)/(\sqrt{\ell}\sigma_q)$. Note that the condition $a>0$
translates to the requirement $\ell< \frac{1}{2E_q}$ as claimed in
the theorem. In case it holds, we have
$$
P\left(x_{\ell}\geq \half \right)\leq \frac{1}{1+\frac{\left( \half
-\ell E_q\right)^2}{\ell \sigma_q^2}}~~,
$$
or put differently,
$$
P(x_\ell < \half)> 1- \frac{1}{1+\frac{\left( \half -\ell
E_q\right)^2}{\ell \sigma_q^2}} = \frac{  \left( \half -\ell
E_q\right)^2}{{\ell \sigma_q^2}+\left( \half -\ell E_q\right)^2}~~,
$$
as stated by the theorem. \hfill $\Box$

\subsection{From Capacity Vector to Coherence}\label{sect_back_to_coherence}

We mentioned earlier that previous work often uses the {\em mutual
coherence} to derive performance bounds on $\ell_1$-reconstructible
supports. The relation between the capacities in $\vect{q}$ and the
inner products between the dictionary atoms,
$|<\vect{d}_i,\vect{d}_j>|$ has been already discussed in \cite{DE1}.
Given a dictionary $\matr{D}$, construct its Gram matrix as $\matr{G}
= \matr{D}^T\matr{D}$. Define the sequence
\begin{eqnarray}
\mu_k=\max_{i\neq k} ~|G_{i,k}| ~~\mbox{for}~~ k\in \Omega.
\end{eqnarray}
Namely, $\mu_k$ is the maximal value on the $k$-th column of
$|\matr{G}|$, disregarding the main diagonal entry. As ~\cite{DE1}
shows, this sequence of values satisfies
\begin{eqnarray}
q_k\leq\frac{\mu_k}{\mu_k+1}.\nonumber
\end{eqnarray}
Thus the condition $ \sum_{k\in\Gamma}q_k<\half$ can be replaced with
$ \sum_{k\in \Gamma}\frac{\mu_k}{\mu_k+1}<\half$, leading of-course,
to weaker bounds. Further relaxation
\begin{eqnarray}
q_k\leq\frac{\mu_k}{\mu_k+1}<\frac{\mu}{\mu+1}
\end{eqnarray}
yields a constant capacity vector with entries of size
$\frac{\mu}{\mu+1}$. Applying Theorem A to this vector we obtain, as
a special case, the classical Theorem \ref{thm_DE1}.

\subsection{Using the Capacity Matrix \matr{Q}}

One problem with the capacity vector $\vect{q}$ is the independence
with which its entries $q_k$ are computed. This implies that one (or
more) of the entries in $\vect{q}$ may become unnecessarily large,
compared to the values obtained in Equation (\ref{C_gamma}), causing
a weaker bound. By working with pairs of such entries, one could in
principle improve the obtained bounds. This leads us to the
following definition:

\begin{defn}
Denote by $\Omega_2$ the set of indices $\Omega_2=\{(i,j)  |~i,j\in
\Omega, i<j\}.$ The upper triangular capacity matrix
$\matr{Q}=\{Q_{i,j}\}$ is the matrix with non-zero elements indexed
by $(i,j)\in\Omega_2$, defined as follows:
$$
Q_{i,j}=\max_{\delta\in Null(\matr{D})}\{\max(\delta_i+\delta_j,
\delta_i-\delta_j)\} \ s.t.\ \ \|\delta\|_1=1.
$$
\end{defn}

Each of these entries can be computed by two independent linear
programming problems of the form
\begin{eqnarray}\nonumber
\left\{\begin{array}{c} \vect{x}_{(i,j)}^+ =
\mbox{Arg}\min_{\vect{x}}~||\vect{x}||_1 ~~\mbox{subject
to}~~\matr{D}\vect{x}=\vect{0} ~\mbox{and}~ x_i +x_j = 1
\\
\vect{x}_{(i,j)}^- = \mbox{Arg}\min_{\vect{x}}~||\vect{x}||_1
~~\mbox{subject to}~~\matr{D}\vect{x}=\vect{0} ~\mbox{and}~ x_i -x_j
= 1
\end{array}\right\}
\end{eqnarray}
 and then assigning $Q_{i,j} =
1/\min(||\hat{\vect{x}}_{(i,j)}^+||_1,||\hat{\vect{x}}_{(i,j)}^-||_1)$.

As in section 3.1, the obtained values $Q_{i,j}$ could be used to
form an improved worst-case bound for Lemma 2.1 and consequently for
Theorem 2.2: Let $\Gamma\subset\Omega$ be a randomly chosen support
of size\footnote{We consider hereafter even support sizes.
Generalization to odd ones is relatively simple, requiring use of one
entry from $\vect{q}$. We omit this discussion for simplicity.}
$\ell=2n$. By definition, the non-zero elements of $\matr{Q}$ satisfy
\begin{eqnarray}\nonumber
\ds \max_{\substack{\delta\in Null(D)\\ \|\delta\|_1=1}}
|\delta_i|+|\delta_j| = Q_{i,j} \le \max_{\substack{\delta\in
Null(D)\\ \|\delta\|_1=1}} |\delta_i| + \max_{\substack{\delta\in
Null(D)\\ \|\delta\|_1=1}} |\delta_j| = q_i+q_j.
\end{eqnarray}
Thus the values $Q_{i,j}$ can be used in the evaluation of an upper
bound on $C_{\Gamma}$. To any partition $\mc{I}$ of $\Gamma$ into
disjoint pairs there corresponds the sum $ \sum_{(k_1,k_2)\in
\mc{I}}Q_{k_1,k_2}$ that bounds the value of $C_{\Gamma}$ from above.
Therefore, $\Gamma$ is $\ell_1$-reconstructible if there exists such
a partition satisfying $\sum_{(k_1,k_2)\in \mc{I}}Q_{k_1,k_2}<
\half$. Naturally, among all such possible partitions, we are
interested in the one that leads to the smallest sum.

Just one glance at the values of $\matr{Q}$ gives a lower bound for
sizes of $\ell_1$-reconstructible subsets: namely, if
$\max(\matr{Q})\leq \frac{1}{\ell}$, then a sum of any $\ell/2$ of
its elements does not exceed $1/2$;  hence any subset of columns of
size up to $\ell$ is guaranteed to be recovered by BP. Conjecture B
below estimates the uncertainty caused by replacing $\max(\matr{Q})$
with $mean(\matr{Q})$. Some numerical techniques based on $\matr{Q}$
are described in section 4.

Here we concentrate again on a theoretical bound that uses
$\matr{Q}$, similar to the one proposed in Theorem A with few
necessary modifications.

We arrange the values $\{Q_{i,j}\;|\;i<j\in \Omega\}$ of the Capacity
matrix in a vector $\matr{Q}^V$. Denote by $E_Q$ the mean value of
$\matr{Q}^V$, and by $\sigma^2_Q$ its variance,
$\sigma^2_Q=\frac{2}{L(L-1)}\sum_{i<j\in \Omega}(Q_{i,j}-E_Q)^2$. The
following statement based on $\matr{Q}$ is similar to the one in
Theorem A:

\begin{thmB}\footnote{This claim is a conjecture since it relies
on a property that is used here without a proof. More on this is
given in Appendix B.}\label{thm_B} For any $1\leq \ell <
\frac{1}{E_Q}$, a support $\Gamma$ of even size $\ell$, sampled
uniformly at random from $\Omega$, admits $\ell_1$-recovery with
probability
\begin{eqnarray}\label{eq:BoundB}
P(\ell) > \frac{ \left( \half
-\frac{\ell}{2}E_Q\right)^2}{\frac{\ell}{2}\sigma_Q^2+\left( \half
-\frac{\ell}{2}E_Q\right)^2}.
\end{eqnarray}

\end{thmB}

Notice that the expression obtained in Equation(\ref{eq:BoundB}) is
the same as the one in (\ref{eq:BoundA}), with $\ell$ replaced by
$\ell/2$. Since $E_Q$ and $\sigma_Q$ refer to pairs, if $E_Q=2E_q$
and $\sigma_Q^2 = 2\sigma_q^2$ the two bounds are the same. However,
as we shall demonstrate in section 5, $E_Q < 2E_q$ and $\sigma_Q^2 <
2\sigma_q^2$ for random dictionaries, implying that this bound is
indeed stronger.

\noindent {\bf Proof:} Fix an even support size $\ell$. In order to
translate the condition  $ \sum_{(i,j)\in \mc{I}}Q_{i,j}<\half$ to a
probabilistic one, we use again the model involving a subset $\Gamma
\subset \Omega$ of size $\ell$ which elements are chosen uniformly
from $\Omega$ without replacement. Also, we let $\mc{I}$ be a random
partition of the index set $\Gamma$ into pairs. Based on these
notions, we define a random variable $x_{\ell}=\sum_{(k_1,k_2)\in
\mc{I}}Q_{k_1,k_2}$. In effect, $x_{\ell}$ is a sum of elements of
$\matr{Q}$ randomly chosen ``without replacement'' in a stronger
sense, i.e. not only the elements are not repeated, but two elements
with common index are not allowed. The probability $P(\ell)$, defined
in the statement of the theorem, is bounded below by $P(x_{\ell}<
\half)$. This bound is not tight, since the support $\Gamma$ is
reconstructible if there exists {\em  some} partition $\mc{I}^{opt}$
such that $\sum_{(k_1,k_2)\in \mc{I}^{opt}}Q_{k_1,k_2}$ drops below
the half, while $P(x_{\ell}< \half)$ is only the probability this
will happen for a {\em random} partition $\mc{I}$.

In order to analyze the variable $x_{\ell}$ we consider a multiset
$\Phi $ of size $\frac{\ell}{2}$ chosen uniformly with replacement
from $\matr{Q}^V$, and define the random variable $y_\ell$ to be its
sum, $y_{\ell}=\sum\Phi.$ Then we have
$\mt{E}(y_{\ell})=\frac{\ell}{2}E_Q$, $var(y_{\ell}) =
\frac{\ell}{2}\sigma_Q^2$.

The expectation of $x_\ell$ equals to that of $y_\ell$, which is
proven in Appendix B. Regarding the variance, we are making an
assumption similar to \ref{X_Y_rel1}:
\begin{equation}\label{X_Y_rel}
var(x_{\ell}) \leq var(y_{\ell}).
\end{equation}

We do not provide its proof and leave it as an open question at this
stage. Empirical verification of this inequality is demonstrated in
Appendix B.

Following the steps of Theorem A, given any real $a>0$, the
one-tailed version of the Tchebychev inequality \cite{Polya} for
$x_{\ell}$ reads
$$
P\left(x_{\ell}\geq \frac{\ell}{2} E_Q+
a\sqrt{\frac{\ell}{2}}\sigma_Q\right)\leq \frac{1}{1+a^2}.
$$

The parameter $a$ is chosen such that $\frac{\ell}{2} E_Q+
a\sqrt{\frac{\ell}{2}}\sigma_Q = \half$, leading to $a=(\half
-\frac{\ell}{2} E_Q)/(\sqrt{\frac{\ell}{2}}\sigma_Q)$, implying that
we should require $\ell< \frac{1}{E_Q}$ to get $a>0$. This leads to
$$
P\left(x_{\ell}\geq \half \right)\leq \frac{1}{1+\frac{\left( \half
-\frac{\ell}{2} E_Q\right)^2}{\frac{\ell}{2} \sigma_Q^2}}~~,
$$
or put differently,
$$
P(x_\ell < \half)> 1-  \frac{1}{1+\frac{\left( \half -\frac{\ell}{2}
E_Q\right)^2}{\frac{\ell}{2} \sigma_Q^2}} = \frac{ \left( \half
-\frac{\ell}{2} E_Q\right)^2}{{\frac{\ell}{2} \sigma_Q^2}+\left(
\half -\frac{\ell}{2} E_Q\right)^2}~~,
$$
as stated in the theorem. \hfill $\Box$


\section{Numerical Algorithms}

Given the capacity vector $\vect{q}$ (or its weaker version as
described in section 3.2) or matrix $\matr{Q}$, we can use Theorems A
and B to predict the $\ell_1$-reconstructible supports, and show
lower bounds of the probability for success as a function of the
support size $\ell$. However, we can alternatively evaluate these
probabilities numerically, provided that there are shortcuts that
avoid the exponential growth in support possibilities. This leads us
to the following two algorithms.

\subsection{A Fast Combinatorial Count Using $\vect{q}$}

Below we propose an algorithm which provides worst-case bounds on
reconstructible support sizes. We would like to establish the
fraction of the total number of supports $\Gamma$ of size $\ell$ that
satisfy $val(C_{\Gamma})< \half$. Testing the sufficient condition $
\sum_{k\in \Gamma}q_k< \half$ for every single $\Gamma$ requires
$\mc{O}(L^\ell)$ flops, which is prohibitive. Instead, we propose to
perform a quantization of the entries of $\vect{q}$ to $d$ distinct
values, and lead to a more reasonable computational process.

Suppose we are given a partition $\Lambda=\{\Lambda_i\}_{i=1}^d$ of
$\Omega$ into $d$ disjoint clusters, such that
$\Omega=\bigcup_{i=1}^{d}\Lambda_i$. The corresponding quantized
values in $\textbf{q}$ are denoted by $ \{q_{\Lambda}^i\}$, each set
to be the maximal in its subset, $ \{q_{\Lambda}^i=\max_{k\in
\Lambda_i}(q_k)\ |\ 1\leq i\leq d\}$.

Given the quantization parameters $\Lambda=\{\Lambda_i,~
q_{\Lambda}^i\}_{i=1}^{d}$, every $\ell$-sized support $\Gamma \in
\Omega$ can be described as the union $\bigcup_{i=1}^{d}\Gamma_i$,
where $\Gamma_i\subseteq \Lambda_i$ is the subset of indices in
$\Gamma$ allocated to the quantized value $q_\Lambda^i$. Thus, the
sum $ \sum_{k\in \Gamma}q_i$ can be replaced by a larger sum,
$\sum_{i=1}^{d}|\Gamma_i|q_\Lambda^i$.

In order to test all possible supports $\Gamma \in \Omega $ of size
$\ell$, a combinatorial count of all sequences $p=(p_1,...,p_d)$ is
performed, such that $0 \leq|p_i|\leq|\Lambda_i|$ and
$\sum_{i=1}^{d}|p_i|=\ell$. For each of these we evaluate
$\sum_{i=1}^{d}|p_i|q^i_{\Lambda}$ and count the relative number of
those\footnote{Each instance must be weighted by the number of its
possible occurrences.} below $\half$. The complexity of such
computation does not exceed $\mc{O}\left((\frac{L}{d})^d\right)$.

As to the choice of the quantization parameters
$\Lambda=\{\Lambda_i,~ q_{\Lambda}^i\}_{i=1}^{d}$, as said above, we
let $q_{\Lambda}^i = \max_{k\in \Lambda_i} q_k$ to guarantee that the
evaluated summations are considering a worst-case scenario. The
clustering is done by an attempt to minimize the function
\begin{eqnarray}
f\left( \{\Lambda_i,~ q_{\Lambda}^i\}_{i=1}^{d} \right) =
\sum_{i=1}^{d}\left(|\Lambda_i|
 q_{\Lambda}^i-\sum_{k\in \Lambda_i}q_k \right).
\end{eqnarray}
The difference $|\Lambda_i| q_{\Lambda}^i-\sum_{k\in \Lambda_i}q_k$
is the quantization error for the elements in the subset
$\Lambda_i$, and the above error simply sums these values.

The minimization of $f\left( \{\Lambda_i,~ q_{\Lambda}^i\}_{i=1}^{d}
\right)$ can be done exhaustively in case $d$ is small -- in our
experiments we have used $d=3$ implying that the above requires
$\mc{O}(L^3)$ flops. For larger values of $d$ a sequential algorithm
that chooses $\Lambda_i$ can be proposed, separating the set $\Omega$
to two parts, and proceeding in a tree and greedy separation scheme.

Computationally, the results of the combinatorial count are very
close to those predicted by Theorem A. Therefore, this method serves
as a supporting evidence for the probabilistic approach taken in
Theorem A, but its numerical output is omitted from our display of
experimental results in section 5.

\subsection{A Sampling Algorithm Using $\matr{Q}$}\label{sect_comp_B}

An alternative to Conjecture B is a direct evaluation of
$\ell_1$-reconstructible supports $\Gamma$ of cardinality $\ell$, by
the following stages:

\begin{itemize}
\item We draw $M\gg L$ such supports $\{\Gamma_i\}_{i=1}^M$.
\item For each $\Gamma_i$ we seek to find a partition $\mc{I}_i$ that
leads to the smallest value of $\sum_{(k,l)\in \mc{I}}Q_{k,l}$. While
finding the best such partition is combinatorial in complexity, we
use an approximate greedy algorithm of complexity $\mc{O}(\ell^2\cdot
log(\ell))$ which computes the following suboptimal partition:
\begin{enumerate}
\item Begin with empty set $\mc{I}$ of pairs.
\item denote by $\matr{Q}_{res}$ the sub-matrix of $\matr{Q}$ which
rows an columns consist of only those indices from $|\Gamma|$ which
do not occur in $\mc{I}$. Retrieve the couple $(i_0 ,j_0), (i_1,
j_1)$ of index pairs which minimize the sum $\matr{Q}(i_0
,j_0)+\matr{Q}(i_1 ,j_1)$ over $\matr{Q}_{res}$.
\item joint the couple $(i_0 ,j_0), (i_1, j_1)$ to $\mc{I}$ and return
to item 2 while $\matr{Q}_{res}$ is nonempty.
\end{enumerate}

Therefore, the algorithm is, in a sense, "second-order greedy", i.e.
at each step the least-sum couple of values from $\matr{Q}$, rather
than least single value, is extracted. Possibly, better algorithms
will improve the performance of this scheme, but we believe it to be
quite close to optimal, while keeping low computational costs. The
fact such partition can be found in $\mc{O}(\ell^2\cdot log(\ell))$
follows from the next combinatorial claim: let $(i^*, j^*)$ be the
index pair of minimal value in submatrix of $\matr{Q}$ supported on
$|\Gamma|$. Then both $i^*, j^*$ necessarily present among indices
$(i_0 ,j_0, i_1, j_1)$ defined above.

\item Given the partition $\mc{I}$, test $\sum_{(k,l)\in
\mc{I}}Q_{k,l} < \half$. Accumulate the relative number of such
occurrences over the collection $\{\Gamma_i\}_{i=1}^M$.
\end{itemize}

\noindent The fact that this method relies on capacity values
implies that the predicted performance is expected to be weaker
compared to the true behavior of BP. Nevertheless, among the various
methods discussed thus far, this method is expected to be the most
optimistic because it uses $\matr{Q}$ and not $\vect{q}$, and also
because it does not build the evaluation through the Tchebychev
inequality that looses also part of the tightness. However, as
opposed to all the other methods described above, this method cannot
claim theoretical correctness of its results.

In the light of similarity of the proposed scheme to the pure
empirical test, we can make a direct comparison of the computational
cost of the two tests. See the details in the Section
\ref{sect_complexity}.


\section{Experimental Results}

\subsection{Test-Cases to Study}

We carry out a number of tests on  each of the three following
dictionaries:
\begin{enumerate}
\item $\matr{D}-Random$ is the  dictionary of size $128\times
256$, which consists of $\ell_2$-normalized random vectors,
independently drawn from the Normal distribution on the unit
sphere. Such a dictionary is often used in numerical experiments
as well as in various applications.

\item $\matr{D}-Spoiled$ is the dictionary $\matr{D}-Random$,
which has undergone an operation designed to create a small set of
columns with high linear dependence. More precisely, we re-generate a
set of $3$ columns as a random linear combination of $12$ other
columns. This dictionary is used to demonstrate the ability of the
{\em capacity-sets} methods to better handle dictionaries with a
non-uniform distribution of inner products.

\item $\matr{D}-DCT$ is the orthonormal pair $[\matr{I},\matr{C}^*]$
of size $128\times 256$, where $\matr{C}$ is the 1-dimensional
Discrete Cosine basis and $\matr{I}$ the identity matrix.
\end{enumerate}

\subsection{Behavior of $\matr{q}$ and  $\vect{Q}$}

As explained earlier, the passage from the capacity vector
$\vect{q}$ to the matrix $\matr{Q}$ was motivated by the fact that
$Q_{i,j}$ provide a lower bound in this context. To exhibit the
numerical behavior of these bounds, we compute the mean and the
variance of the family of ratios
\begin{eqnarray}
R_{k,l}=\frac{Q_{k,l}}{q_k+q_l}\ ~~\mbox{for}~~k\neq l\in \Omega .
\end{eqnarray}
The mean and variance of these ratios for the three test cases is
given in Table \ref{TableI}.

\begin{table}[htbp]
\begin{center}
\scriptsize{
\begin{tabular}{||c||c|c||}
\hline Dictionary & $\mt{E}\left(R\right)$ &
$\sigma\left(R\right)$
\\\hline \hline $\matr{D}-Random$ & 0.7175 & 0.0008 \\ \hline
$\matr{D}-Spoiled$ & 0.7154  & 0.001 \\ \hline
$\matr{D}-DCT$ & 0.6509 & 0.0109 \\
\hline
\end{tabular}}
\end{center}
\caption{Behavior of the {\em capacity-sets} $\vect{q}$ and
$\matr{Q}$ by evaluating the mean and variance of the
ratios.\label{TableI}}
\end{table}

As these figures show, we earn up to $30\%$ of the upper bound value
by upgrading to Capacity Matrix from the Capacity Vector. This ratio
between the two bounds for the corresponding indices is very stable,
as seen from the low values of the standard deviation
$\sigma\left(R\right)$.

To display the power of Conjecture B, we show that $E_Q < 2E_q$ and
either $\sigma_Q^2 < 2\sigma_q^2$ or $\sigma_Q^2 \ll E_Q^2$. The
corresponding values for various dictionaries are presented in the
table below.

\begin{table}[htbp]
\begin{center}
\scriptsize{\hspace{-3in}
\begin{tabular}{||c||c|c||c|c||}
\hline Dictionary & $E_Q$ & $2E_q$ & $\sigma_Q^2$  & $2\sigma_q^2$ \\
\hline
\hline $\matr{D}$-Random $32\times  128$ & 0.2329  & 0.3179 & 0.5849e-3 & 0.8252e-3\\
\hline $\matr{D}$-Random $64\times  128$ & 0.1695  & 0.2345 & 0.1405e-3 & 0.1654e-3\\
\hline $\matr{D}$-Random $128\times 256$ & 0.1235  & 0.1721 & 0.4511e-4 & 0.5652e-4\\
\hline $\matr{D}$-DCT $64\times     128$ & 0.1687  & 0.2586 & 0.4732e-3 & 0.0112e-3\\
\hline $\matr{D}$-DCT $128\times    256$ & 0.1265  & 0.1943 & 0.4070e-3 & 0.4144e-5\\
\hline
\end{tabular}}
\end{center}
\caption{Comparison of mean and variance of capacity
sets.\label{TableII}}
\end{table}

Notice that for the $\matr{D} - DCT$ dictionary the variance of the
capacity vector is smaller than that of the Capacity matrix, due to
the special structure of this dictionary. Nevertheless, as seen later
in the results section, Conjecture B predicts BP success on support
sizes larger than those allowed by Theorem A.


\subsection{Compared Methods}

We perform a number of computations, applying various methods for
the estimation of BP performance on the given dictionaries. The
results are expressed via a set of Estimation Functions,
$EF:\Omega\ra \mt{R}$, which value at $\ell\in \Omega$ is the
predicted percentage of $\ell$-sized supports which admit recovery
by $\ell_1$-norm optimization. The EFs considered are the
following:
\begin{enumerate}
\item EF-emp - The standard empirical test on the dictionary. This
test is done by drawing $1,000$ random supports for each cardinality
$\ell$, generating a corresponding signal, and solving the BP per
each. EF-emp is obtained by showing the relative number of successes
in recovering the support.

\item EF-CB - the classical coherence-based upper bound
$\half(1+\frac{1}{\mu})$, provided by the Theorem \ref{thm_DE1}.

\item EF-thmA - expresses the results of the Theorem A, EF-thmA
$(\ell)= P(\ell)$ as defined in the statement of the theorem. The
values are computed from $\vect{q}$ of the dictionary.

\item EF-thmB - expresses the results of the Conjecture B, computed
from the capacity matrix $\matr{Q}$ of the dictionary.

\item EF-compB - The results of the sampling algorithm based on
$\matr{Q}$, which results support the estimation of Conjecture B (see
section 4.2).

\item EF-GB - The Grassmanian upper bound, computed by the formula
for the Classical Bound using the ideal coherence
$\mu=\sqrt{\frac{L-N}{N(L-1)}}$.
\end{enumerate}

\noindent This last EF deserves more explanation: Among all possible
dictionaries of size $N\times L$, the Grasssmanian frame is the one
leading to the smallest possible coherence $\mu =
\sqrt{\frac{L-N}{N(L-1)}}$ \cite{SH}. Thus, this leads to the most
optimistic worst-case bound. When the dictionary is ``un-balanced'',
implying a large spread of inner-products in the Gram-matrix, we
know that the {\em mutual-coherence}-bound deteriorates
dramatically. Thus, by using the Grassmanian Bound, we test what is
the best achievable coherence-based performance behavior for the
same dictionary size.

\subsection{Complexity Analysis of the Methods}\label{sect_complexity}

We argue the usefulness of Capacity-based numerical algorithms for
an evaluation of a given dictionary $\matr{D}$. To that end, we
consider the computational complexity of each method listed in
previous section.

\begin{enumerate}
\item EF-emp - The standard empirical test of $\matr{D}$ is conveyed as
follows: for each support size $\ell$, pick $M>>L$ random subsets
$\Gamma$ of columns of size $\ell$. For each $\Gamma$, generate a
signal with  random coefficients vector supported on $\Gamma$ and
test if BP will recover the support. Since in practice maximal
relevant size $\ell$ is proportional to $L$, the computational
complexity of this test is $\mc{O}(M\cdot L\cdot C_{LP}(L))$, where
$C_{LP}(L)$ denotes the complexity of linear programming algorithm
for problem of size $L$.

\item EF-CB requires the computation of $\mu$, which takes
$\mc{O}(L\cdot N)$ flops.

\item EF-thmA - To employ results of the Theorem A, the capacity vector
$\matr{q}$ is computed in ( $\mc{O}(L\cdot C_{LP}(L))$), and then for
each $\ell$ the probability $P(\ell)$, defined in the statement of
Theorem A, is computed in $\mc{O}(L)$. Overall complexity -
$\mc{O}(L^2+L\cdot C_{LP}(L))=\mc{O}(L\cdot C_{LP}(L))$.

\item EF-thmB - To employ results of the Conjecture B, the capacity vector
$\matr{q}$ is computed in ( $\mc{O}(L^2\cdot C_{LP}(L))$), and then
for each $\ell$ the probability $P(\ell)$, defined in the statement
of Conjecture B, is computed in $\mc{O}(L^2)$. Overall complexity -
$\mc{O}(L^3+L^2\cdot C_{LP}(L))=\mc{O}(L^2\cdot C_{LP}(L))$.

\item EF-compB - Our heaviest (and best-performance) algorithm conducts
a semi-empirical test: for each support size $\ell$, pick $M>>L$
random subsets of columns of size $\ell$, and employ the analysis
detailed in \ref{sect_comp_B}. The computational cost of single
support treatment is $\mc{O}(\ell^2\cdot log(\ell))$. Overall
complexity is $\mc{O}(L^2\cdot C_{LP}(L) +M\cdot L^2\cdot log(\L))$.
\end{enumerate}

\noindent As seen from the analysis above, only the EF-compB has
non-negligible computational complexity. When comparing EF-emp and
EF-compB, we can concentrate on the relative complexities of linear
programming solver versus the $\mc{O}(\ell^2\cdot log(\ell))$ of the
partition algorithm, and the benefit of the later is evident.

\subsection{Comparison Results}

Figure 1 presents the obtained graphs of the various EF-s functions
described above, for the three dictionaries described at the top of
this section. As we see from the left-side graphs in the figures, for
all the dictionaries the empirically established support size which
admits BP recovery is at least $40$ columns. Note that this relative
number of columns is also predicted in \cite{DT1}, however, this
holds true only asymptotically (for dictionaries of growing sizes)
and for specific random dictionaries.

Returning to statements which hold for our modest size of $128\times
256$, we notice that the estimation made by the sampling algorithm
based on the Capacity Matrix (EF-compB) is much better than the
Classical bound, established so far in the literature. The difference
is especially high for the D-Spoiled dictionary, which reflects the
fact that methods based on {\em capacity sets} manage well the
non-uniform distribution of inner products.
\begin{figure}[htbp]
    \begin{tabular} {ll}
        \hspace{-1.6in}\includegraphics[width=0.55\textwidth]
        {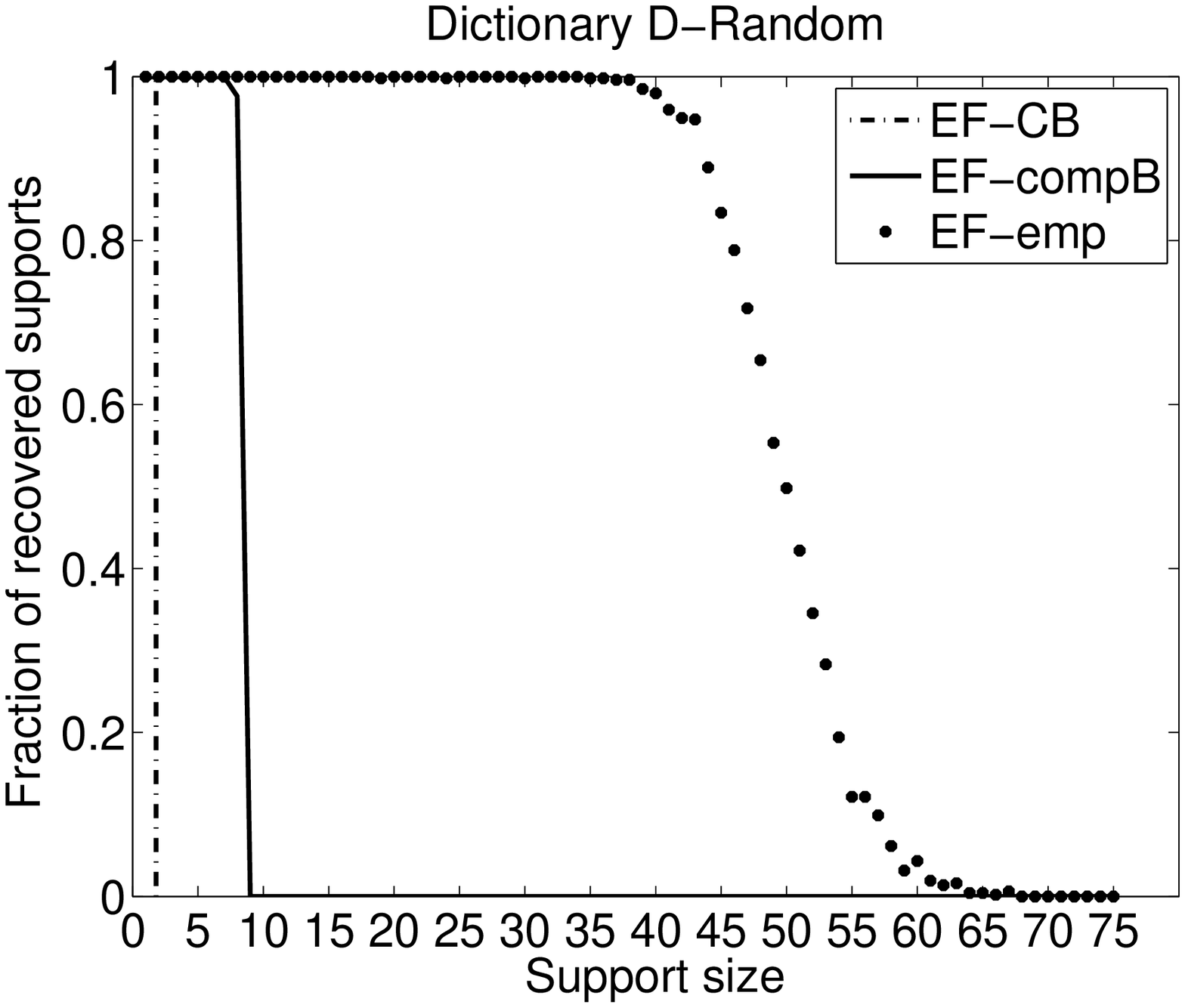} &
        \hspace{-0in}\includegraphics[width=0.55\textwidth]
        {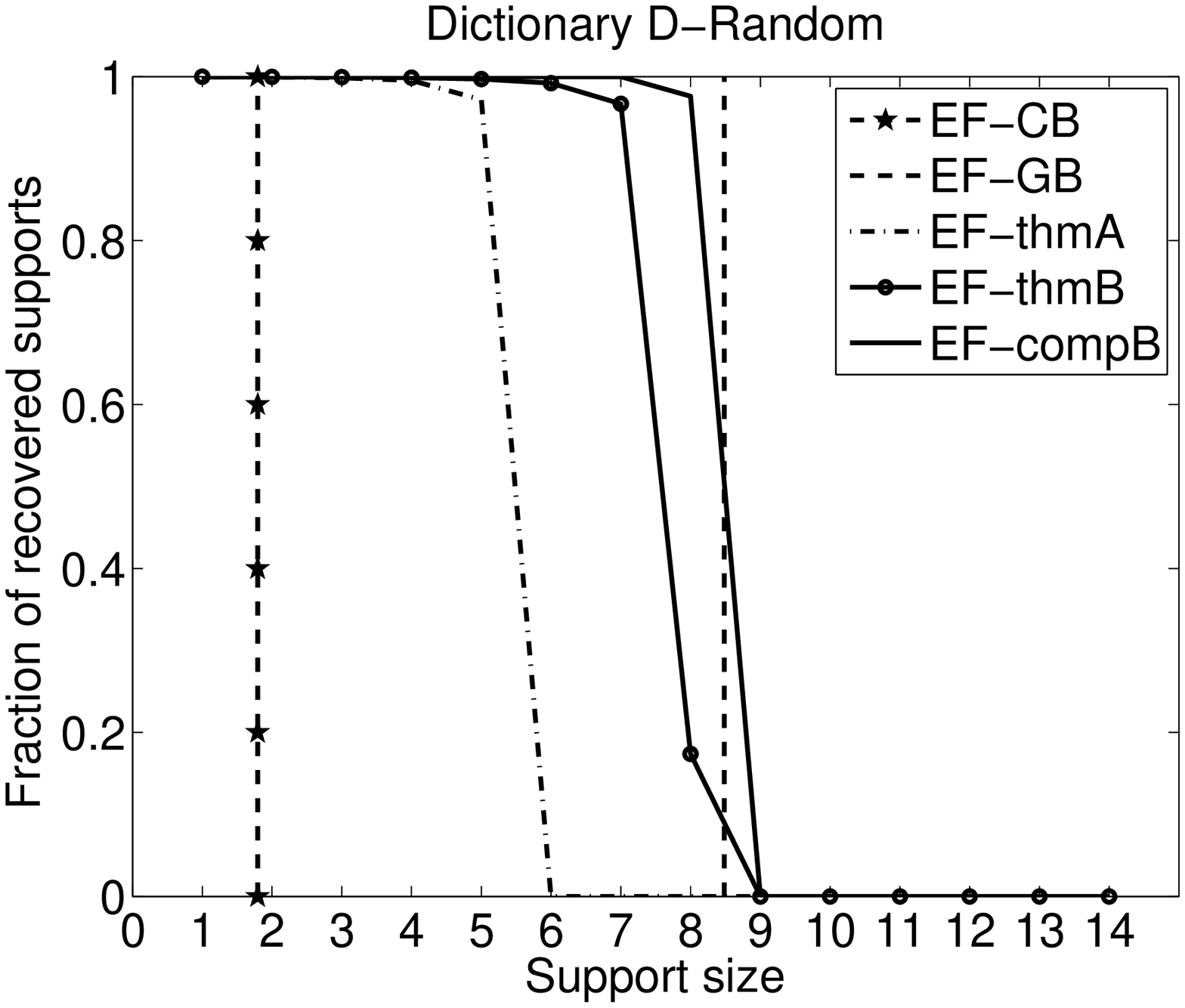} \\

       \hspace{-1.6in}\includegraphics[width=0.55\textwidth]
        {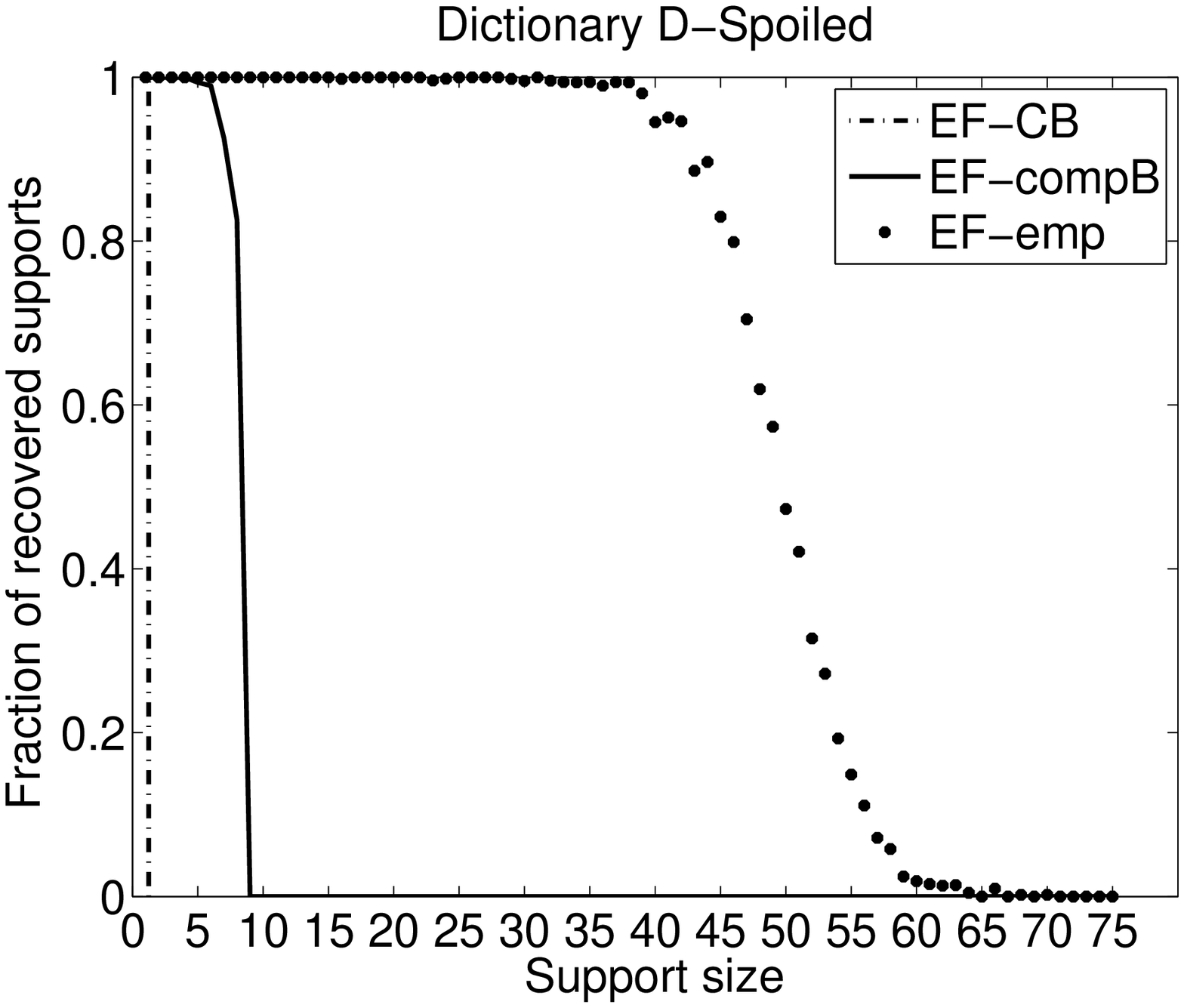} &
       \hspace{-0in}\includegraphics[width=0.55\textwidth]
       {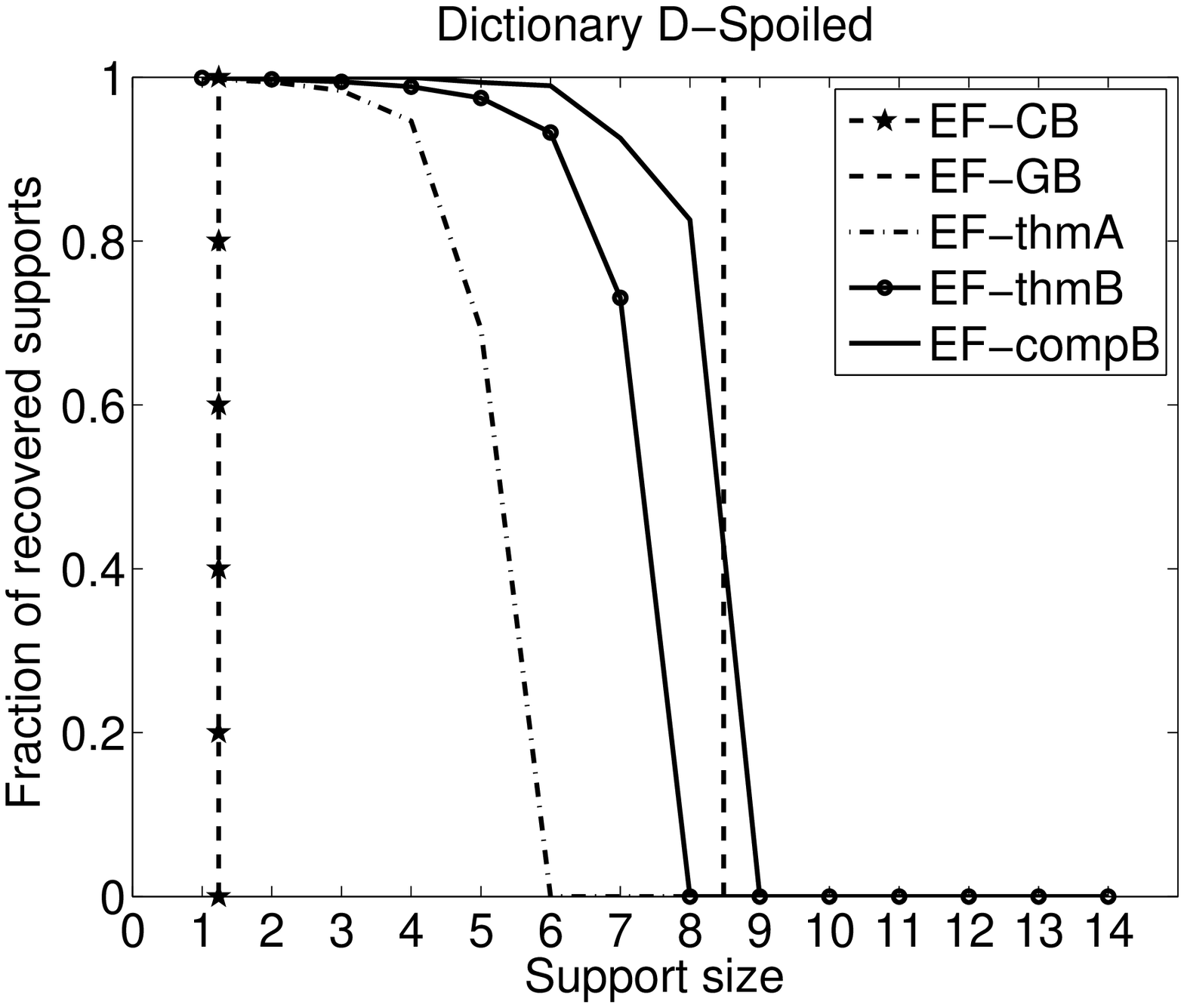} \\

      \hspace{-1.6in}\includegraphics[width=0.55\textwidth]
        {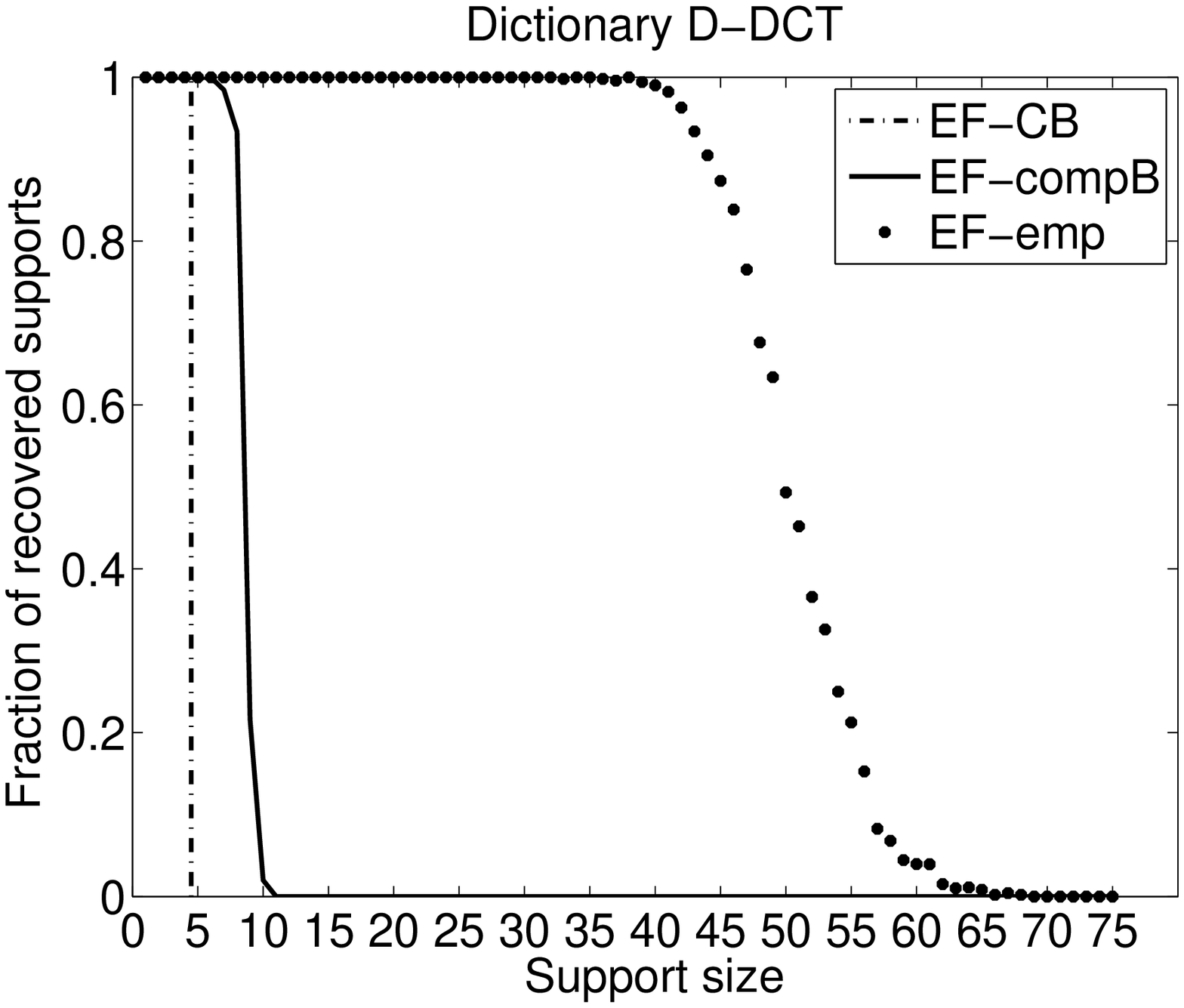} &
       \hspace{-0in}\includegraphics[width=0.55\textwidth]
       {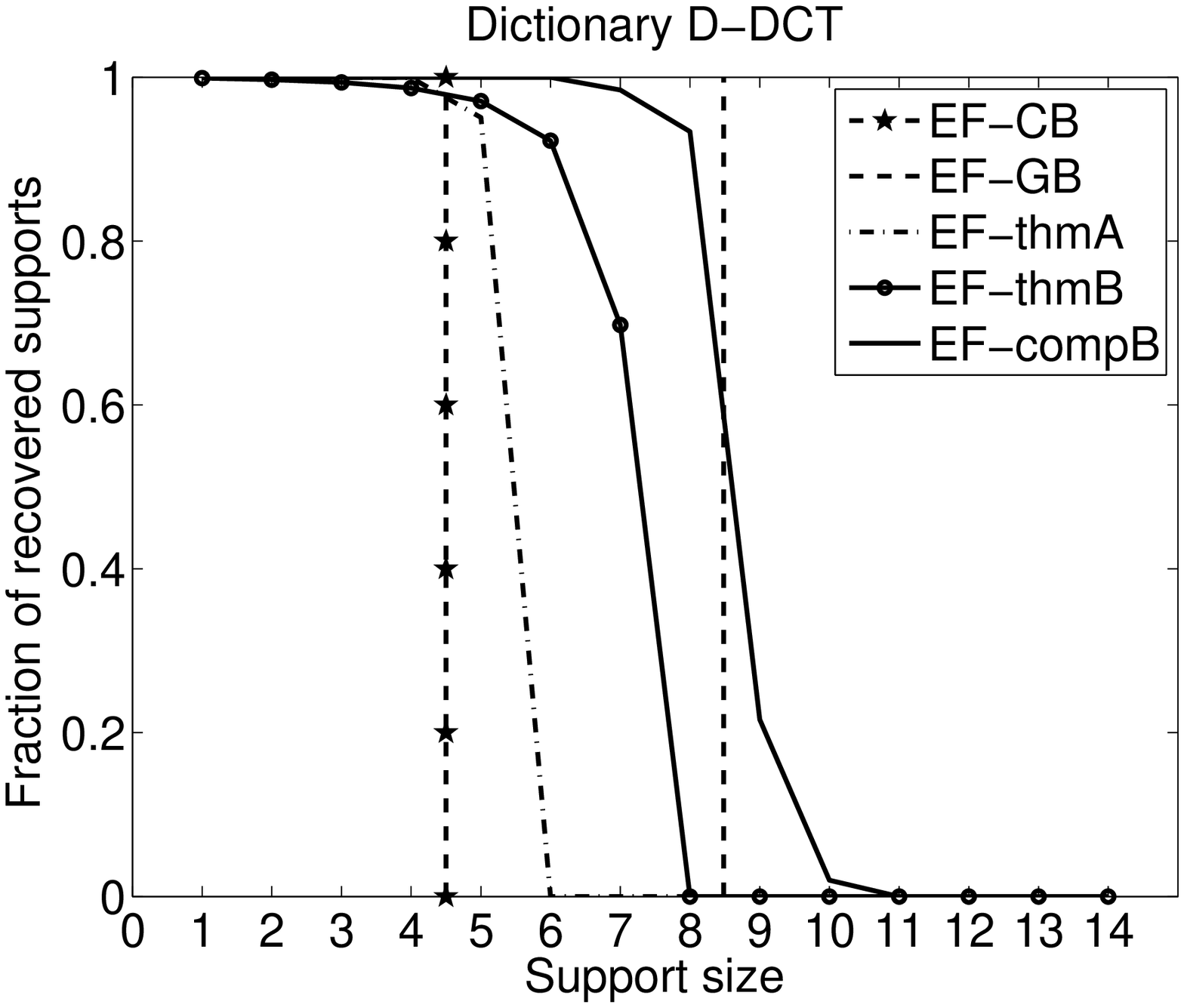} \\

        \end{tabular}

   \caption{\scriptsize{Estimation Functions for various dictionaries
            of size $128\times 256$.}}
\end{figure}

On the right side of each figure we display various method developed
in this work. Noticeably, the results of Conjecture B(EF-thmB) are
stronger than those of Theorem A (EF-thmA), which is explained by the
benefit of using the Capacity Matrix rather than the Capacity Vector.
This benefit is expressed in the ratio values given in Tables
\ref{TableI}, \ref{TableII} and explained thereafter. Apparently,
Conjecture B does not express the full power of the Capacity Matrix
estimation, since the sampling algorithm based on its values
(EF-compB) outperforms EF-thmB by $15-20\%$. This algorithm produces
values which are quite close to the Grassmanian Bound, the best
possible bound one can hope to obtain using coherence-based
estimation for the given dictionary size. We do not have enough
information to explain the fact that values of EF-compB and of
Grassmanian bound nearly coincide for all the dictionaries discussed
here (and additional ones examined during the work); Discovering the
reason underlying this connection may be a lead to important insights
regarding the Basis Pursuit performance.


\section*{Appendix A}
\renewcommand{\theequation}{A-\arabic{equation}}
\setcounter{equation}{0}

\noindent We prove the claim \ref{X_Y_rel1}.
\begin{thmC}
For the two random variables, $x_{\ell}$ and $y_{\ell}$, defined in
\ref{defn_xl_yl}, the following relations between the first and
second moments hold:
\begin{equation}
\mt{E}(x_{\ell})=\mt{E}(y_{\ell})~~\mbox{and}~~ var(x_{\ell})\leq
var(y_{\ell}).
\end{equation}
\end{thmC}

\noindent {\bf Proof:} We begin by introducing some notation. Fix the
support size $1\leq \ell \leq L$. For any $1\leq k\leq \ell$, we
denote by $\mc{C}^k_{\ell}$ the collection of all $\ell$-sized
non-ordered multisets of indices from $\Omega$ (with repetitions),
which have precisely $k$ distinct elements each. For instance,
$\{1,4,5,4,7\}$ and $\{5,1,7,4,4\}$ are two distinct elements of
$\mc{C}^4_{5}$. Such multiset will be sometimes referred to as "index
set". Also, we define $\mc{D}^n_{\ell}=\mc{C}^\ell_{\ell}\cup
\mc{C}^{\ell-1}_{\ell}\cup...\cup \mc{C}^{\ell-n}_{\ell} $, the
collection of all $\ell$-sized multisets having at least $\ell-n$
distinct elements.

In this notation, $x_{\ell}$ is a random variable with uniform
distribution over the domain $\mc{D}^0_{\ell}$, which admits value
$\sum_{k\in \Lambda}q_k$ on a given element $ \Lambda\in
\mc{D}^0_{\ell}$. The variable $y_{\ell}$ has the same definition on
a larger domain $\mc{D}^{\ell-1}_{\ell}$, containing the domain of
$x_{\ell}$. Therefore, we treat both $x_{\ell}$ and $y_{\ell}$ as
restrictions of the same uniformly distributed random variable $x$ on
the corresponding domains: $x_{\ell}=x_{|\mc{D}^{0}_{\ell}}$,
$y_{\ell}=x_{|\mc{D}^{\ell-1}_{\ell}}$. In the proof we use the
following basic property of the variance:

\begin{prop}\label{var_union}
Let $z$ be a random variable defined over a domain given as the
disjoint union $\mc{D}=\mc{D}_1\cup \mc{D}_2\cup...\cup \mc{D}_n$,
with uniform distribution. Denote $v=var(z_{|\mc{D}}),
v_i=var(z_{|\mc{D}_i}),s_i=|\mc{D}_i|$. Then $\ds
v=\frac{\sum_{i=1}^n s_i v_i}{\sum_{i=1}^n s_i}$.
\end{prop}

\noindent \textbf{Part 1.} The expectation of the random variable
$x$ restricted to $\mc{D}^0_{\ell}$ is computed by
$$\begin{array}{l}\ds
\mt{E}(x_{|\mc{D}^0_{\ell}})=\frac{1}{|\mc{D}^0_{\ell}|}\sum_{\Lambda\in
\mc{D}^0_{\ell}}\sum_{k\in \Lambda}q_k.
\end{array}$$
This sum contains $|\mc{D}^0_{\ell}|\cdot\ell$ elements, and for each
$j\in \Omega$, $q_j$ appears in it the same number of times.
Therefore, each $q_j$ appears $|\mc{D}^0_{\ell}|\frac{\ell}{L}$
times, and we have $\ds
\mt{E}(x_{|\mc{D}^0_{\ell}})=\frac{\ell}{L}\sum_{k\in \Omega}q_k=\ell
E_q$. The mean of $x_{|\mc{D}^{\ell-1}_{\ell}}$ is computed
similarly:
$$\begin{array}{l}\ds
\mt{E}(x_{|\mc{D}^{\ell-1}_{\ell}})=\frac{1}{|\mc{D}^{\ell-1}_{\ell}|}
\sum_{\Lambda\in\mc{D}^{\ell-1}_{\ell}}\sum_{k\in \Lambda}q_k.
\end{array}$$
Here each $q_j$ appears $|\mc{D}^{\ell-1}_{\ell}|\frac{\ell}{L}$
times, and we have $\ds
\mt{E}(x_{|\mc{D}^{\ell-1}_{\ell}})=\frac{\ell}{L}\sum_{k\in
\Omega}q_k=\ell E_q$.

This proves our first claim, $\mt{E}(x_{\ell})=\mt{E}(y_{\ell})$.
For the rest of the proof, where only the variance of the two
variables is considered, we assume w.l.g. that the expectation of
$x_{\ell}$ and $y_{\ell}$ is zero (in the light of equality
$var(z)=var(z-\mt{E}(z)$ for any random variable $z$), that is
$E_q=0$.

\noindent \textbf{Part 2.} We consider the extension of $x$, defined
so far on domain comprising of distinct $\ell$-sized index sets, to
the domain where each such set may appear any finite number of times.
$x$ still has a uniform distribution over this collection. Thus, a
disjoint union of two or more (non-necessarily distinct) index sets
is a sub-domain to which $x$ may be restricted.

For any $0\leq n < \ell$, we define two disjoint unions
$$\begin{array}{l} \ds
\mc{A}_n=\bigcup_{\Gamma\in \mc{D}^n_{\ell-1}}\{\Gamma\cup\{j\}\ |\ j\in \Gamma\},
\\ \ds \mc{B}_n=\bigcup_{\Gamma\in
\mc{D}^n_{\ell-1}}\{\Gamma\cup\{j\}\ |\ j\in \Omega\}
\end{array}$$
(In the definition of $\mc{A}_n$, the set $\Gamma\cup\{j\}$ is added
to the collection one time for each appearance of $j$ in $\Gamma$.)

Let  $\Lambda\in\mc{C}^k_{\ell}$ be a set which contains distinct
indices $j_1,...,j_k$ with multiplicities $m_1,..,m_k$ (so that
$\sum_{i=1}^k m_i=\ell$). For each $1\leq i\leq k$, $\Lambda$ is
obtained in $\mc{A}_n$ $m_i-1$ times in the form $\Gamma\cup\{j_i\}$
for an appropriate $\Gamma=\Gamma_i\in \mc{C}^k_{\ell-1}$ (this claim
also holds vacuously for $m_i=1$). Therefore, the number of copies of
$\Lambda$ in $\mc{A}_n$ equals $\sum_{i=1}^k(m_i-1)=\ell-k$. Also,
$\Lambda$ appears in $\mc{B}_n$ precisely once for each
$j_1,...,j_k$, in the form $\Gamma\cup\{j_i\}$ (for an appropriate
$\Gamma=\Gamma_i$ each time). Therefore, $\mc{B}_n$ contains $k$
copies of  $\Lambda$.

Denote a disjoint union of $a$ distinct copies of some collection $\mc{C}$ by
$a\cdot\mc{C}$. Then we can write $\mc{A}_n,\mc{B}_n$ as
\begin{equation}\label{union_A}
\mc{A}_n=0\cdot \mc{C}^{\ell}_{\ell}\cup 1\cdot
\mc{C}^{\ell-1}_{\ell}\cup...\cup n\cdot \mc{C}_{\ell}^{\ell-n}
\end{equation}
\begin{equation}\label{union_B}
\mc{B}_n=\ell\cdot \mc{C}^{\ell}_{\ell}\cup (\ell-1)\cdot
\mc{C}^{\ell-1}_{\ell}\cup...\cup (\ell-n)\cdot \mc{C}_{\ell}^{\ell-n}
\end{equation}
We prove the following inequality:
$$ var(x_{|\mc{B}_n})\leq var(x_{|\mc{A}_n}).$$
Since $E_q=0$ by our assumption, the expectations of
$x_{|\mc{A}_n}$ and $x_{|\mc{B}_n}$ also equal zero: by the
argument similar to one presented in the first part of the proof,
$\mt{E}(x_{|\mc{A}_n})=\mt{E}(x_{|\mc{B}_n})=\ell\cdot E_q$. Thus
we have

$$\begin{array}{l} \ds
var(x_{|\mc{A}_n})= \frac{1}{|\mc{D}^n_{\ell-1}|}\sum_{\Gamma\in
\mc{D}^n_{\ell-1}}\frac{1}{\ell-1}\sum_{j\in\Gamma}(\sum_{k\in\Gamma}q_k+q_j)^2.
\end{array}$$
For the brevity of the argument we introduce the notation
$\ds q_{\Gamma}=\sum_{k\in\Gamma}q_k$.\\
Then $var(x_{|\mc{A}_n})$ reads as
\begin{eqnarray}var(x_{|\mc{A}_n}) & = & \nonumber
\frac{1}{|\mc{D}^n_{\ell-1}|}\sum_{\Gamma\in
\mc{D}^n_{\ell-1}}\frac{1}{\ell-1}\sum_{j\in\Gamma}(q^2_{\Gamma}+q^2_j+
2q_{\Gamma}q_j)=
\\ \nonumber & =&\frac{1}{|\mc{D}^n_{\ell-1}|}\sum_{\Gamma\in
\mc{D}^n_{\ell-1}}q^2_{\Gamma}+\frac{1}{\ell-1}\sum_{j\in\Gamma}(q^2_j+
2q_{\Gamma}q_j).
\end{eqnarray}

Similarly, we have
\begin{eqnarray}
var(x_{|\mc{B}_n}) \nonumber & = & \frac{1}{|\mc{D}^n_{\ell-1}|}
\sum_{\Gamma\in
\mc{D}^n_{\ell-1}}\frac{1}{L}\sum_{j\in\Omega}(\sum_{k\in\Gamma}q_k+q_j)^2=
\\ \nonumber & = &
\frac{1}{|\mc{D}^n_{\ell-1}|}\sum_{\Gamma\in
\mc{D}^n_{\ell-1}}q^2_{\Gamma}+\frac{1}{L}\sum_{j\in\Omega}(q^2_j+
2q_{\Gamma}q_j).
\end{eqnarray}

The summand $\ds \frac{1}{|\mc{D}^n_{\ell-1}|}\sum_{\Gamma\in
\mc{D}^n_{\ell-1}}q^2_{\Gamma}$ appears in both expressions hence
cancels out. We consider the term $\ds
\frac{1}{|\mc{D}^n_{\ell-1}|}\sum_{\Gamma\in \mc{D}^n_{\ell-1}}
\frac{1}{\ell-1}\sum_{j\in\Gamma}q^2_j$ in $var(x_{|\mc{A}_n})$.
The element $q^2_a$ appears in it same number of times for every
$a\in \Omega$. Hence $\ds
\frac{1}{|\mc{D}^n_{\ell-1}|}\sum_{\Gamma\in \mc{D}^n_{\ell-1}}
\frac{1}{\ell-1}\sum_{j\in\Gamma}q^2_j=\frac{1}{L}\sum_{a\in\Omega}q^2_a$.
By same argument, in  the expression of  $var(x_{|\mc{B}_n})$ we
have $\ds \frac{1}{|\mc{D}^n_{\ell-1}|}\sum_{\Gamma\in
\mc{D}^n_{\ell-1}}
\frac{1}{L}\sum_{j\in\Omega}q^2_j=\frac{1}{L}\sum_{a\in\Omega}q^2_a$,
hence this quadratic term also cancels out. In the light of these
observations, we obtain
$$\begin{array}{l} \ds
var(x_{|\mc{A}_n})-var(x_{|\mc{B}_n})=\frac{2}{|\mc{D}^n_{\ell-1}|}
\sum_{\Gamma\in \mc{D}^n_{\ell-1}}
q_{\Gamma}(\frac{1}{\ell-1}\sum_{i\in\Gamma}q_i-\frac{1}{L}
\sum_{j\in\Omega}q_j).
\end{array}$$
Here we substitute again $q_{\Gamma}$ for $\ds\sum_{i\in\Gamma}q_i$
and recall $\ds \frac{1}{L}\sum_{j\in\Omega}q_j=E_q=0$. Thus, we
have
$$\ds var(x_{|\mc{A}_n})-var(x_{|\mc{B}_n})=
\frac{2}{(\ell-1)|\mc{D}^n_{\ell-1}|}\sum_{\Gamma\in
\mc{D}^n_{\ell-1}}q^2_{\Gamma}\geq 0.
$$

In order to use this result for the proof of the theorem, we make the
following observations : Denote $v_n=var(x_{|\mc{C}^n_{\ell}})$ and
$s_n=|\mc{C}^n_{\ell}|$. By virtue of the decomposition
(\ref{union_A}), $var(x_{|\mc{A}_n})$ can be written as $\ds
var(x_{|\mc{A}_n})=\frac{\sum_{i=0}^n i\cdot s_{\ell-i}v_{\ell-i}}{
\sum_{i=0}^n i\cdot s_{\ell-i}}$ (see Proposition \ref{var_union}).
Similarly, we have $\ds var(x_{|\mc{B}_n})=\frac{\sum_{i=0}^n
(\ell-i)\cdot s_{\ell-i}v_{\ell-i}}{ \sum_{i=0}^n (\ell-i)\cdot
s_{\ell-i}} $. We compute the coefficients of $v_i$ in the expression
\begin{equation}\nonumber
var(x_{|\mc{A}_n})-var(x_{|\mc{B}_n}) =\frac{\sum_{i=0}^n i\cdot
s_{\ell-i}v_{\ell-i}}{ \sum_{i=1}^n i\cdot s_{\ell-i}}
-\frac{\sum_{i=0}^n (\ell-i)\cdot s_{\ell-i}v_{\ell-i}}{ \sum_{i=1}^n
(\ell-i)\cdot s_{\ell-i}}.
\end{equation}
For any $0\leq k\leq n$, the coefficient of $v_{l-k}$ is
$$\begin{array}{l} \ds
\dfrac{1}{Den}s_{\ell-k}\left( k\sum_{i=1}^n (\ell-i)\cdot
s_{\ell-i}-(\ell-k)\sum_{i=1}^n i\cdot s_{\ell-i}\right) \\
\nonumber \hspace{2in} = \dfrac{1}{Den}\ell\cdot
s_{\ell-k}\sum_{i=0}^n (k-i)s_{\ell-i},
\end{array}$$
with
$$\begin{array}{l} \ds
Den=\sum_{i=1}^n i\cdot s_{\ell-i}\cdot \sum_{i=1}^n (\ell-i)\cdot s_{\ell-i} .
\end{array}$$

 We denote $\ds \alpha_{\ell-k}=\ell\sum_{i=0}^n (k-i)s_{\ell-i}$, for
$1\leq k\leq n$, in order to write the above difference as
\begin{equation}
0\leq var(x_{|\mc{A}_n})-var(x_{|\mc{B}_n})=\dfrac{1}{Den}\sum_{k=0}^n
\alpha_{\ell-k}s_{\ell-k} v_{\ell-k}.
\end{equation}
The  constant $\dfrac{1}{Den}$ is positive, since $n< \ell$. Thus,it
can be omitted while preserving the inequality:
\begin{equation}\label{var_diff_2}
0\leq \sum_{k=0}^n \alpha_{\ell-k}s_{\ell-k} v_{\ell-k}.
\end{equation}
The coefficients in this expression have the two following properties:
\begin{enumerate}
\item $\sum_{k=0}^n s_{\ell-k}\alpha_{\ell-k}=0.$
\item $\forall j, \alpha_{j-1}-\alpha_{j}=\ell\sum_{i=0}^ns_{\ell-i}$.
\end{enumerate}
To show the first equality, we consider the sum in (1) as the
linear combination of the elements $s_{\ell-i}s_{\ell-j}$,
$i,j=0,...,n$. The coefficient of $s_{\ell-i}s_{\ell-i}$ is zero
for any $i$. For any $i\neq j$,  $s_{\ell-i}s_{\ell-j}$ appears
just in two components of the sum above, namely,
$s_{\ell-i}\alpha_{\ell-i}$ and $s_{\ell-j}\alpha_{\ell-j}$.
Specifically, $\alpha_{\ell-i}$ contains the summand
$\ell(i-j)s_{\ell-j}$, and $\alpha_{\ell-j}$ contains the summand
$\ell(j-i)s_{\ell-i}$, therefore in the sum
$s_{\ell-i}\alpha_{\ell-i}+s_{\ell-j}\alpha_{\ell-j}$ the
coefficient of $s_{\ell-i}s_{\ell-j}$ is zero. The second property
follows from the definition of $\alpha_{i}$. In the light of the
first property, \ref{var_diff_2} can be written as
\begin{equation}\label{var_diff_3}
(\sum_{k=1}^n \alpha_{\ell-k}s_{\ell-k})v_{\ell}\leq
   \sum_{k=1}^n \alpha_{\ell-k}s_{\ell-k} v_{\ell-k}.
\end{equation}

Equipped with these observations, we prove, by induction on $n$,
the inequality
$$
var(x_{|\mc{D}^0_{\ell}})\leq var(x_{|\mc{D}^n_{\ell}}).
$$
for any $n=1,...,\ell-1$. The theorem follows for $n=\ell-1$. By
Proposition \ref{var_union}, $\ds
var(x_{|\mc{D}^{n}_{\ell}})= \frac{\sum_{i=0}^n
s_{\ell-i}v_{l-i}}{\sum_{i=0}^n s_{\ell-i}}$, and
$var(x_{|\mc{D}^0_{\ell}})$ is just $v_{\ell}.$ Thus we need
to prove
$$
v_{\ell}\leq\frac{\sum_{i=0}^n s_{\ell-i}v_{l-i}}{\sum_{i=0}^n
s_{\ell-i}},
$$
or
\begin{equation}\label{var_claim}
(\sum_{i=1}^n s_{\ell-i})v_{\ell}\leq\sum_{i=1}^n s_{\ell-i}v_{l-i}.
\end{equation}
For $n=1$, \ref{var_diff_3} reads as
$$
\alpha_{\ell-1}s_{\ell-1}v_{\ell}\leq
\alpha_{\ell-1}s_{\ell-1}v_{\ell-1}.
$$ Here $\alpha_{\ell-1}=\ell s_{\ell}>0$, thus we obtain
the inequality
$$
s_{\ell-1}v_{\ell}\leq s_{\ell-1}v_{\ell-1},
$$
as required. Now, we assume by induction that inequality
\ref{var_claim} holds up to $n-1$ and prove for $n$. We use
(\ref{var_diff_3}):
$$
(E1):\ \ \ (\sum_{k=1}^n \alpha_{\ell-k}s_{\ell-k})v_{\ell}\leq
\sum_{k=1}^n \alpha_{\ell-k}s_{\ell-k} v_{\ell-k}.
$$
This inequality undergoes a series of transformations designed to
bring it to the form of \ref{var_claim}.

First, we have $\alpha_{\ell-1}<\alpha_{\ell-2}$. Since
$v_{\ell}\leq v_{\ell-1}$ by the proof for $n=1$, we have an
inequality
$$
(d1):\ \ \ (\alpha_{\ell-2}-\alpha_{\ell-1})s_{\ell-1}v_{\ell}\leq
(\alpha_{\ell-2}-\alpha_{\ell-1})s_{\ell-1}v_{\ell-1}
$$
Adding $(d1)$ to the inequality $(E1),$ we arrive at
\begin{eqnarray} (E2):\ \  \
& & \left(\alpha_{\ell-2}(s_{\ell-1}+s_{\ell-2})+ \sum_{k=3}^n
\alpha_{\ell-k}s_{\ell-k}\right)v_{\ell}\leq \nonumber \\
\nonumber & & \hspace{0.3in} \leq
\alpha_{\ell-2}(s_{\ell-1}v_{\ell-1}+s_{\ell-2}v_{\ell-2})+
\sum_{k=3}^n \alpha_{\ell-k}s_{\ell-k} v_{\ell-k}.
\end{eqnarray}

Second, by induction assumption for $n=2$ we have the inequality
$$(s_{\ell-1}+s_{\ell-2})v_{\ell}\leq s_{\ell-1}v_{\ell-1}+
s_{\ell-2}v_{\ell-2}.$$
Also, $\alpha_{\ell-2}\leq \alpha_{\ell-3}$ as noticed earlier.
Then we can construct the next inequality in order to add it to
$(E2)$:
$$
(d1):\ \ \
(\alpha_{\ell-3}-\alpha_{\ell-2})(s_{\ell-1}+s_{\ell-2})v_{\ell}\leq
(\alpha_{\ell-3}-\alpha_{\ell-2})(s_{\ell-1}v_{\ell-1}+s_{\ell-2}v_{\ell-2})
$$
This results in the following expression:
\begin{eqnarray}(E3):\ \ \ & &
\left(\alpha_{\ell-3}\sum_{i=1}^3 s_{\ell-i}+
\sum_{k=4}^n \alpha_{\ell-k}s_{\ell-k}\right)v_{\ell}\leq \nonumber
\\ \nonumber & & \hspace{0.3in} \leq \alpha_{\ell-3}\sum_{i=1}^3
(s_{\ell-i}v_{\ell-i})+ \sum_{k=4}^n \alpha_{\ell-k}s_{\ell-k}
v_{\ell-k}.
\end{eqnarray}
In this fashion we make $n-1$ steps resulting in the inequality
$$
(E(n)):\ \ \ (\alpha_{\ell-n}\sum_{i=1}^n s_{\ell-i})v_{\ell}\leq
\alpha_{\ell-n}\sum_{i=1}^n s_{\ell-i}v_{\ell-i}
$$
Notice that $\alpha_{\ell-n}$ is positive:
$\alpha_{\ell-n}=s_{\ell-n}\ell(ns_{\ell}+(n-1)s_{\ell-1}+...
+s_{\ell-n+1})$.
Thus, we obtain the desired result. As mentioned, the theorem
follows for $n=\ell-1$. \hfill $\Box$

\section*{Appendix B}
\renewcommand{\theequation}{B-\arabic{equation}}
\setcounter{equation}{0}

We prove the equality of expectations
\begin{equation}
\mt{E}(x_{\ell}) = \mt{E}(y_{\ell}),
\end{equation}
for random variables $x_{\ell}$ and $y_{\ell}$ defined in the proof
of Conjecture B. Recall that $y_\ell$ is a sum of $\frac{\ell}{2}$
values from $\matr{Q}$, uniformly distributed over this matrix,
therefore $\mt{E}(y_{\ell})=\frac{\ell}{2}\mt{E}_Q.$ We show
$\mt{E}(x_{\ell})=\frac{\ell}{2}\mt{E}_Q$, too, by considerations of
symmetry, similar to those used in the proof of Theorem A, part 1.

Namely, we consider a totality $\mc{P}_\ell$ of partitions of all
$\ell$-sized supports $\Lambda\subset \Omega$, into ordered pairs of
indices. An element in this collection is therefore a pair
$(\Lambda,\mc{I}_{\Lambda})$. We clarify that the index sets
$\Lambda\subset \Omega$ are chosen without repetitions and up to a
permutation of their elements. Now, let $(i,j)$ be an ordered pair of
indices from $\Omega$. We argue that the number of appearances of
this pair in the elements of $\mc{P}_\ell$ does not depend on choice
of $i$ and $j$. Indeed, this number is just the size of the
collection $\mc{P}_{\ell-2}$, built for submatrix of $\matr{Q}$ with
$i$-th and $j$-th rows and columns missing.

Since $x_\ell(\Lambda,\mc{I}_{\Lambda})$ is the sum $\sum_{(i,j)\in
\mc{I}_{\Lambda}}\matr{Q}(i,j)$, we conclude that all the elements
$\matr{Q}(i,j)$ contribute to the value of $x_\ell$ with equal
probability, hence $\mt{E}(x_{\ell})=\frac{\ell}{2}\mt{E}_Q$ as
desired.

Now we provide an empirical evidence to the claim
\begin{equation}\label{App_B_var}
 var(x_{\ell}) \leq var(y_{\ell})
\end{equation}
Statistical data below supports this inequality. While the variance
of $y_{\ell}$ is known precisely, for $x_{\ell}$ we estimate it by
drawing $10^4$ random subsets of indices for each support size up to
half the signal dimension of the dictionary. Results are presented in
Figure 2. The computation is carried out for a number of dictionary
sizes on dictionary $\matr{D}$-Random. As can be seen from these
figures, the gap between $var(x_{\ell})$ and $var(y_{\ell})$ is
roughly proportional to the support size.

Same experiments on dictionary $\matr{D}$-DCT display different
results: the variance of both variables coincides. As number of
samples grows, we observe that the difference of variance values, for
all support sizes, tends to zero. We conclude that for this specific
dictionary, \ref{App_B_var} is an equality.

\begin{figure}[htbp]
    \begin{tabular} {ll}
        \hspace{-1.6in}\includegraphics[width=0.5\textwidth]
        {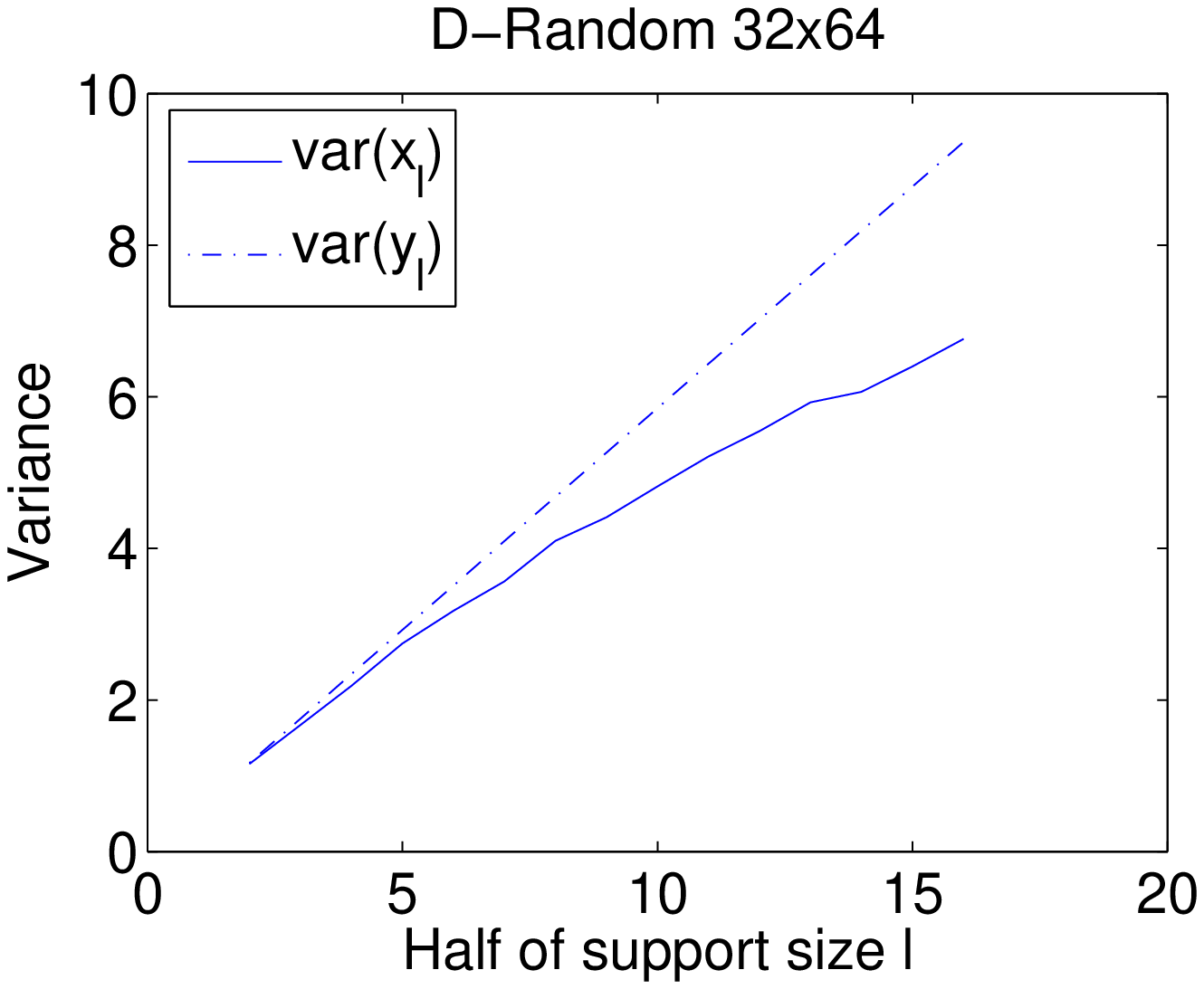} &
        \hspace{-0in}\includegraphics[width=0.5\textwidth]
        {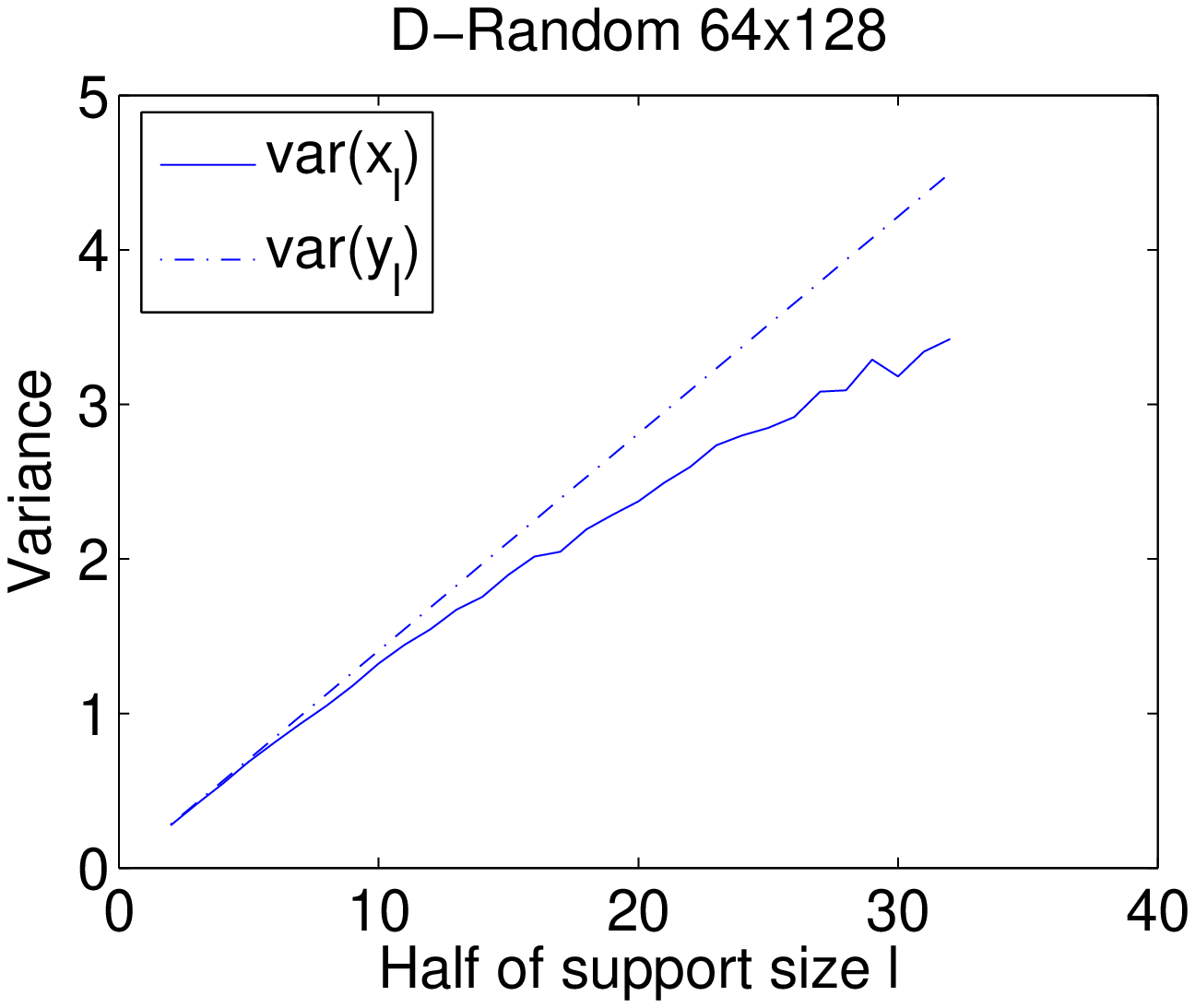}\\
       \hspace{-1.6in}\includegraphics[width=0.5\textwidth]
        {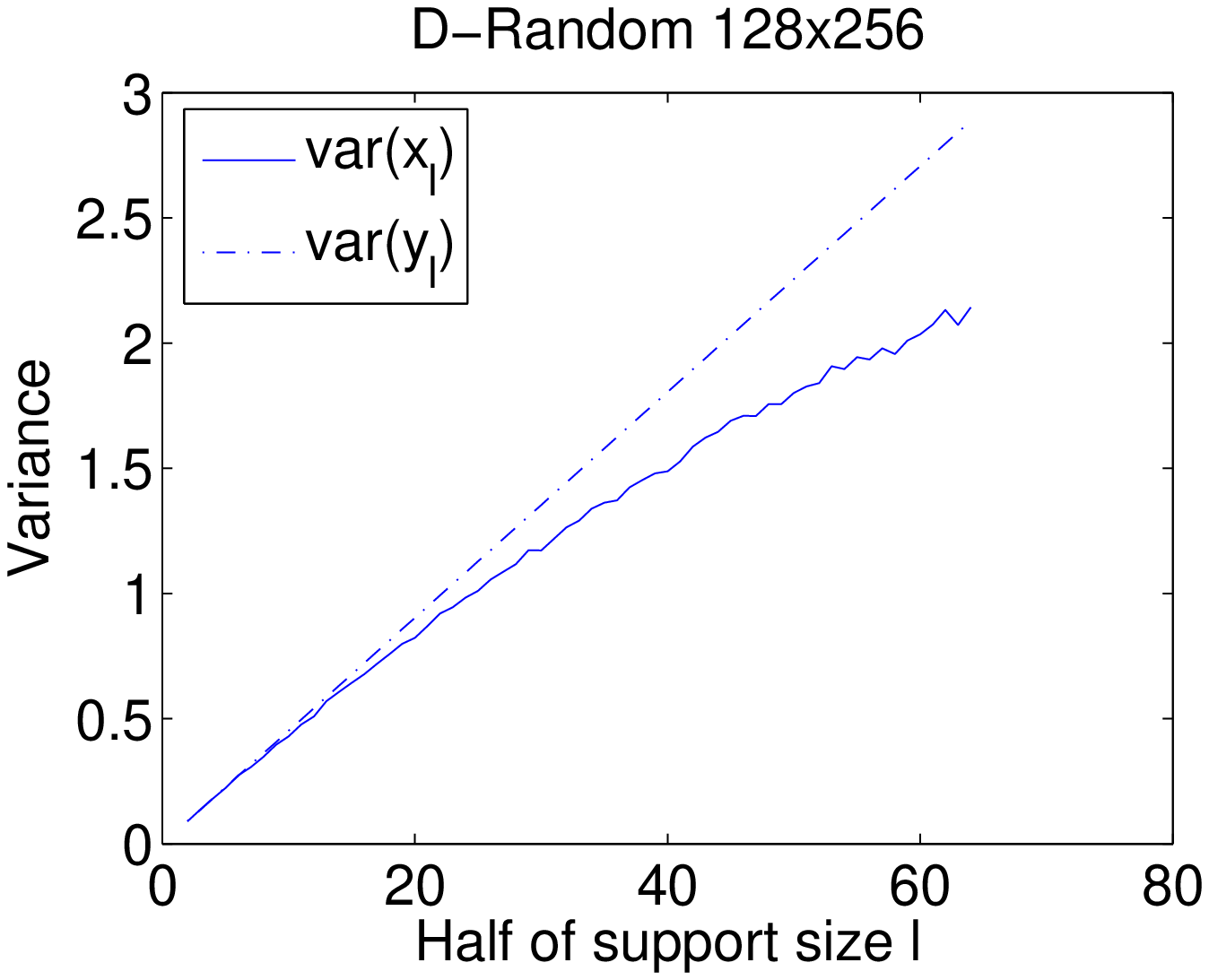} \\
        \end{tabular}
    \label{fig_Graph_variance2}
   \caption{\scriptsize{ The variances of $x_{\ell}$ and $y_{\ell}$
            (scaled by $10^{3}$)}}
\end{figure}


\end{document}